\documentclass[a4paper,11pt]{amsart}

\usepackage{amsmath}
\usepackage{amssymb}
\usepackage{amsthm}
\usepackage{mathtools}
\usepackage[dvipdfmx]{graphicx}
\usepackage{cite} 
\usepackage{color}
\usepackage{tikz-cd} 
\usepackage{comment}
\usepackage{blkarray}
\usepackage{todonotes}
\usepackage{microtype}
\usepackage{bm}
\usepackage{hyperref}
\usepackage{enumitem}
\setlist[enumerate,1]{label={\upshape (\arabic*)}}

 \newtheorem{theorem}{Theorem}[section]
 \newtheorem{lemma}[theorem]{Lemma}
 \newtheorem{proposition}[theorem]{Proposition}

\newtheorem{corollary}[theorem]{Corollary}

\theoremstyle{definition}
\newtheorem{definition}[theorem]{Definition}
\newtheorem{remark}[theorem]{Remark}
\newtheorem{example}[theorem]{Example}

\numberwithin{equation}{section}

\newcommand{\NN}{\mathbb{N}}
\newcommand{\MCG}{\mathrm{MCG}}
\newcommand{\SC}{\Sigma^\circ}
\newcommand{\GC}{G^\circ}
\newcommand{\fC}{f^\circ}
\newcommand{\ZC}{Z^\circ}

\newcommand{\fg}{\mathfrak{g}}
\newcommand{\marginal}[1]{\marginpar{\tiny #1}}

\title[Polynomial conjugacy invariants of pseudo-Anosov maps]
{Polynomial conjugacy invariants of pseudo-Anosov maps 
derived from train track representatives}

\author[E. Kin]{%
    Eiko Kin
}
\address{%
      Center for Education in Liberal Arts and Sciences, The University of Osaka, Toyonaka, Osaka 560-0043, Japan
}
\email{%
        kin.eiko.celas@osaka-u.ac.jp
}

\author[Y. Mizutani]{%
    Yuta Mizutani
}
\address{%
      Department of Mathematics, Graduate School of Science, The University of Osaka, Toyonaka, Osaka 560-0043, Japan
}
\email{%
        yuta.mizutani.math@gmail.com
}



\date{
November 7, 2026
}

\begin{document}

\begin{abstract}
We investigate 
two kinds of polynomials associated with pseudo-Anosov maps.
Each kind consists of three polynomials that are conjugacy invariants of pseudo-Anosov mapping classes 
and these polynomials contain more information than 
pseudo-Anosov stretch factors. 
The first kind  is defined by the train track representatives.
The second kind is obtained from the first kind by incorporating additional data that encode the transitions  
among the unpunctured singularities of the invariant foliations. 
The Birman-Brinkmann-Kawamuro polynomials coincide 
with those of the second kind. 
As an application, we prove that 
the Birman-Brinkmann-Kawamuro polynomials 
are conjugacy invariants of pseudo-Anosov mapping classes 
which is also proved by Dodge-Kawamuro, independently. 
\end{abstract}

\maketitle

\section{Introduction}
\label{section_introduction}

\begin{figure}[htbp]
\begin{center}
\includegraphics[width=3.8in]{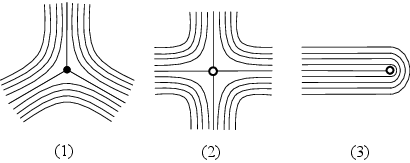}
\caption{
The unstable/stable foliations near singularities.  
($\bullet$ and $\circ$ denote unpunctured and punctured singularities, 
respectively.)  
(1) unpunctured $3$-pronged, 
(2) punctured $4$-pronged and 
(3) punctured $1$-pronged singularities.}
\label{fig_foliation}
\end{center}
\end{figure}

Let $\Sigma= \Sigma_{g,s}$ be an orientable surface 
of genus $g$ with $s$ punctures.    
The {\it mapping class group} $\mathrm{MCG}(\Sigma)$  is 
the group of isotopy classes of 
orientation preserving self-homeomorphisms on $\Sigma$ 
preserving punctures set-wise. 
Assume that $3g-3+s \ge 1$. 
The Nielsen-Thurston classification of surface homeomorphisms  
states that 
any homeomorphism on  $\Sigma$ 
is isotopic to one of the three types, 
periodic, reducible, and pseudo-Anosov \cite{Thurston88}. 
Here, a homeomorphism 
$f: \Sigma \rightarrow \Sigma$  is {\it periodic} if 
some power of $f$ is the identity map on $\Sigma$, 
and {\it reducible} if 
there exists a collection of mutually and non-homotopic essential 
simple closed curves $C_1, \dots, C_n$ in $\Sigma$, 
called {\it reducing curves} for $f$, 
such that 
$f(C_1 \cup \dots \cup C_n)= C_1 \cup \dots \cup C_n$. 
A homeomorphism $f: \Sigma \rightarrow \Sigma $ is {\it pseudo-Anosov}  
if  there exist a constant $\lambda >1$ 
and a pair of transverse measured foliations 
$(\mathcal{F}^u, \mu^u)$ and 
$(\mathcal{F}^s, \mu^s)$ 
such that 
$$f(\mathcal{F}^u, \mu^u) = (\mathcal{F}^u, \lambda \mu^u) \  \hspace{2mm}\mbox{and}\ \hspace{2mm}
f(\mathcal{F}^s, \mu^s) = (\mathcal{F}^s, \lambda^{-1} \mu^s).$$
We call  $\lambda$  the {\it stretch factor} of 
the pseudo-Anosov mapping class $[f]$.  
The constant $\lambda$ is an algebraic integer 
that is a conjugacy invariant of $[f]$ in $\mathrm{MCG}(\Sigma)$ 
(see \cite{FLP79}), 
and $\mathcal{F}^u$ and $\mathcal{F}^s$ are called 
{\it unstable} and {\it stable foliations}. 
They are singular foliations in $\Sigma$ and 
admit a common finite set of singularities  
that include punctures of $\Sigma$ (Figure \ref{fig_foliation}). 
The set of singularities consists of  periodic points of $f$. 
We define a {\it singularity orbit} of $f$ 
as the periodic orbit of a singularity.

Let $f: \Sigma \rightarrow \Sigma$ be a pseudo-Anosov map 
and $\mathcal{S}$  the set of all unpunctured singularities. 
We say that $f$ is \textit{fully punctured} 
if $\mathcal{S}= \emptyset$, that is, 
the set of singularities equals the set of punctures of $\Sigma$. 
We have a fully punctured pseudo-Anosov map 
$f^{\circ}: \Sigma^{\circ} \rightarrow \Sigma^{\circ}$ 
 obtained from $f: \Sigma \rightarrow \Sigma$ 
 by removing all unpunctured singularities, that is, 
$f^{\circ}:= f|_{\Sigma^{\circ}}: 
\Sigma^{\circ} \rightarrow \Sigma^{\circ}$ 
where $\Sigma^{\circ}:= \Sigma \setminus \mathcal{S}$. 

In this paper, we investigate 
two kinds of polynomials associated with pseudo-Anosov maps.
Each kind consists of three polynomials that are conjugacy invariants of pseudo-Anosov mapping classes. 
See the polynomials (1)(2)(3) in Table \ref{table_two-kinds}. 
To introduce the first kind, 
recall that for the computation of 
the stretch factor of $[f]$, 
one uses a train track $\tau$ 
which is a finite $C^1$ graph embedded in $\Sigma$. 
If $\tau$ is {\it $f$-invariant}, then $f$ determines 
a graph map $\mathfrak{t}: \tau \rightarrow \tau$ 
that takes vertices to vertices, edges to edge paths. 
See Section \ref{subsection_weight-spaces}. 
This graph map $\mathfrak{t}: \tau \rightarrow \tau$ 
gives a  combinatorial description of the pseudo-Anosov map and 
is called a {\it train track representative} of $[f]$. 
Let $\mathcal{E}(\tau)$ be the set of all edges of $\tau$, and 
let $|\mathcal{E}(\tau)|$ be its cardinality.
The graph map $\mathfrak{t}$  naturally induces a linear map 
$$\mathfrak{t}_*: {\Bbb R}^{|\mathcal{E}(\tau)|} \rightarrow 
{\Bbb R}^{|\mathcal{E}(\tau)|}$$
using a correspondence between an $i$th edge of $\tau$ 
and a standard unit vector of ${\Bbb R}^{|\mathcal{E}(\tau)|}$ 
whose $i$th component is $1$. 
Then the largest real root of the characteristic polynomial 
of  $\mathfrak{t}_*$ on ${\Bbb R}^{|\mathcal{E}(\tau)|}$ 
is the stretch factor of $[f]$. 
Note that this polynomial is not a conjugacy invariant of $[f]$, 
since it depends on the choice of the train track. 
Instead of the vector space  ${\Bbb R}^{|\mathcal{E}(\tau)|}$, 
we consider the {\it weight space} $W(\tau)$. 
This is  a subspace of 
${\Bbb R}^{|\mathcal{E}(\tau)|}$ 
which consists of vectors that satisfy the switch conditions. 
Then $\mathfrak{t}_*(W(\tau)) \subset W(\tau)$ 
and we have  the restriction map  
$\mathfrak{t}_*: W(\tau) \rightarrow W(\tau)$. 
The first kind of polynomials  consists of 
the following three characteristic polynomials given in Theorem 
\ref{theorem_first-kind}.

\begin{theorem}
\label{theorem_first-kind}
Let $f: \Sigma \rightarrow \Sigma$ be a pseudo-Anosov map. 
The characteristic  polynomial  $\chi(\mathfrak{t}_*)$ 
of $\mathfrak{t}_*: W(\tau) \rightarrow W(\tau)$ 
factors as 
$\chi(\mathfrak{t}_*)= \chi(\mathfrak{t}_*|_{Z(\tau)}) \cdot 
\chi(\overline{\mathfrak{t}_*})$. 
Polynomials $\chi(\mathfrak{t}_*)$, $\chi(\mathfrak{t}_*|_{Z(\tau)})$ 
and $\chi(\overline{\mathfrak{t}_*})$ are
conjugacy invariants of $[f] \in \mathrm{MCG}(\Sigma)$.  
\end{theorem}

We call the three polynomials in Theorem \ref{theorem_first-kind} 
the {\it train track polynomials} of $[f]$. 
In Section \ref{subsection_symplectic-form}, 
we define the Thurston symplectic form $\omega$ on $W(\tau)$ 
which is not a symplectic form in general. 
We denote by $Z(\tau)$,  the radical of $\omega$, 
and consider a direct sum decomposition 
$W(\tau) \simeq Z(\tau) \oplus W(\tau)/Z(\tau)$. 
The map $\mathfrak{t}_*: W(\tau) \rightarrow W(\tau)$ 
preserves the form $\omega$, and 
the quotient space 
$W(\tau)/Z(\tau)$ inherits the form $\omega$. 
Hence, we have the restriction map 
$\mathfrak{t}_*|_{Z(\tau)}: Z(\tau) \rightarrow Z(\tau)$ 
and the induced map 
$\overline{\mathfrak{t}_*}: W(\tau)/Z(\tau) \rightarrow W(\tau)/Z(\tau)$. 
The characteristic polynomials of 
$\mathfrak{t}_*|_{Z(\tau)}$ and $\overline{\mathfrak{t}_*}$ 
occur  as factors of $\chi(\mathfrak{t}_*)$ in 
Theorem \ref{theorem_first-kind}. 

In the proof of Theorem \ref{theorem_first-kind}, 
we need the notion of {\it unstable/stable stars} introduced by 
Hironaka-Tsang \cite{HironakaTsang25}. 
We also need an explicit description of the radical of the 
Thurston symplectic form given by \cite{HironakaTsang25}.

\begin{table}[hbtp]
\caption{
Two kinds of polynomial conjugacy invariants: 
Train track polynomials and BBK polynomials.}
\label{table_two-kinds}
\begin{center}
\begin{tabular}{|l|l|}
\hline
\hspace{5mm}  {\bf  Train track polynomials}  & 
\hspace{10mm}{\bf  BBK polynomials} 
\\
\hline 
 {\small Train track representative}  
$\mathfrak{t}: \tau \rightarrow \tau$ & 
{\small Train track map}
$\mathfrak{g}: G \rightarrow G$ 
\\
\hline 
$\mathfrak{t}_*: W(\tau) \rightarrow W(\tau)$ & 
$\mathfrak{g}_*: W(G) \rightarrow W(G) $
\\
\hline
{\small (1) $\chi(\mathfrak{t}_*)$ of 
$\mathfrak{t}_*: W(\tau) \rightarrow W(\tau)$} 
& 
{\small (1) $ \chi(\mathfrak{g}_*)$ of 
$\mathfrak{g}_*: W(G) \rightarrow W(G)$}
\\
\hline
{\small (2) Radical polynomial $\chi(\mathfrak{t}_*|_{Z(\tau)})$} & 
{\small (2) Puncture polynomial $\chi(\mathfrak{g}_*|_{Z(G)})$}
\\
\hline
\multicolumn{2}{|c|}
{{\small (3) Symplectic polynomial $\chi(\overline{\mathfrak{t}_*}) 
= \chi(\overline{\mathfrak{g}_*})$}} 
\\ 
\hline
{\small $\chi(\mathfrak{t}_*)= 
\chi(\mathfrak{t}_*|_{Z(\tau)}) \chi(\overline{\mathfrak{t}_*})$} & 
{\small $\chi(\mathfrak{g}_*)= 
\chi(\mathfrak{g}_*|_{Z(G)}) \chi(\overline{\mathfrak{g}_*})$}
\\
\hline
\end{tabular}
\end{center}
\end{table}

The second kind of polynomials is based on the Bestvina-Handel 
algorithm \cite{BestvinaHandel95,Brinkmann00}. 
See Section \ref{section_BBKpolynomials}. 
Let $f: \Sigma \rightarrow \Sigma$ be a homeomorphism 
on a punctured surface $ \Sigma= \Sigma_{g,s}$ for $s \ge 1$. 
Then one can find an {\it induced graph map}  
$\mathfrak{g}': G' \rightarrow G'$ of $[f]$, 
where $G'$ is embedded in $\Sigma$ 
and is homotopy equivalent to $\Sigma$. 
From the graph map $\mathfrak{g}': G' \rightarrow G'$, 
the algorithm determines the Nielsen-Thurston type of $f$. 
In case $[f]$ is pseudo-Anosov, 
the algorithm finds an {\it efficient graph map} 
$\mathfrak{g}: G \rightarrow G$ that is an induced graph map of $[f]$, 
where $G$ is embedded in  $\Sigma$ 
and is homotopy equivalent to $\Sigma$. 
We call the resulting graph map 
$\mathfrak{g}: G \rightarrow G$ 
a {\it train track map} of $[f]$. 
In \cite{BestvinaHandel95} 
Bestvina-Handel construct a train track $\tau= \tau(G)$ 
together with a surjection $\pi: \tau \rightarrow G$. 
They also construct 
a train track representative $\mathfrak{t}: \tau \rightarrow \tau$ 
that is an extension of 
$\mathfrak{g}: G \rightarrow G$. 
Then $\tau= \tau(G)$ has a property such that 
for each {\it even type} (resp. {\it odd type}) vertex $v$ of $G$, 
$\pi^{-1}(v)$ forms an infinitesimal even-gon (resp. odd-gon) 
corresponding to an unpunctured 
even-pronged (resp. odd-pronged) singularity 
of $\mathcal{F}^u$ and $\mathcal{F}^s$. 
See Definition \ref{definition_vertex-type}.

Let $\mathfrak{g}_*: {\Bbb R}^{|\mathcal{E}(G)|} 
\rightarrow {\Bbb R}^{|\mathcal{E}(G)|} $ 
be a linear map 
induced by the train track map $\mathfrak{g}: G \rightarrow G$. 
Then $\mathfrak{g}_*$ restricts to a subspace 
$W(G):= \pi_*(W(\tau))$, 
where $\pi_{*}: {\Bbb R}^{|\mathcal{E}(\tau)|} 
\rightarrow {\Bbb R}^{|\mathcal{E}(G)|}$ 
is a surjection induced by $\pi: \tau \rightarrow G$. 
Denoting  
$\pi_*|_{W(\tau)}: W(\tau) \rightarrow W(G)$ 
by $\pi_*: W(\tau) \rightarrow W(G)$, 
we have the  exact sequence 
$$0 \longrightarrow \ker \pi_* \longrightarrow W(\tau) 
\stackrel{\pi_*}{\longrightarrow} W(G) \longrightarrow 0$$

The Thurston symplectic form $\omega$ on $W(\tau)$ gives rise to 
a skew-symmetric form $\overline{\omega}$ on $W(G)$. 
Birman-Brinkmann-Kawamuro \cite{BBK12} studied 
 $\mathfrak{g}_*: W(G) \rightarrow W(G)$, 
the restriction map 
$\mathfrak{g}_*|_{Z(G)}: Z(G) \rightarrow Z(G)$ 
and the induced map 
$\overline{\mathfrak{g}_*}: W(G)/Z(G) \rightarrow W(G)/Z(G)$ 
that are linear maps parallel to those defining the train track polynomials. 
We call their characteristic polynomials the {\it BBK polynomials} 
that are of the second kind in this paper.

In Section \ref{subsection_unpunctured}, 
we prove that the dimension of the subspace $\ker \pi_*$ 
is equal to the number of even type vertices of $G$. 
The map 
$\mathfrak{t}_*: W(\tau) \rightarrow W(\tau)$ 
restricts to  $\ker \pi_*$. 
The restriction  $\mathfrak{t}_*|_{\ker \pi_*}: \ker \pi_* 
\rightarrow \ker \pi_*$ 
encodes the transitions among unpunctured even-pronged 
singularities of the pseudo-Anosov map. 
We prove in Theorem \ref{theorem_kerpi*_coninv} that 
the characteristic polynomial 
$\chi(\mathfrak{t}_*|_{\ker \pi_*})$ 
is a conjugacy invariant of 
$[f] \in \mathrm{MCG}(\Sigma)$. 

The following result states a relation between 
the two kinds of polynomials.

\begin{theorem}
\label{theorem_characterization}
Let $f: \Sigma \rightarrow \Sigma$ be a pseudo-Anosov map, 
where $\Sigma= \Sigma_{g,s}$ with $s \ge 1$. 
Let $\mathfrak{g}\colon G\to G$ be a train track map and $\mathfrak{t}\colon \tau \to \tau$ the train track representative of $[f]$. 
Then we have the following. 
\begin{enumerate}
\item
$\chi(\mathfrak{t}_*)
=\chi(\mathfrak{g}_*)\chi(\mathfrak{t}_*|_{\ker\pi_*})$. 

\item 
$\chi(\mathfrak{t}_*|_{Z(\tau)})=\chi(\mathfrak{g}_*|_{Z(G)})
\chi(\mathfrak{t}_*|_{\ker\pi_*})$. 

\item 
$\chi(\overline{\mathfrak{t}_*})= \chi(\overline{\mathfrak{g}_*})$
\end{enumerate}
In particular, if $f$ is a fully punctured pseudo-Anosov map, 
then 
$\chi(\mathfrak{t}_*)=\chi(\mathfrak{g}_*)$, 
$\chi(\mathfrak{t}_*|_{Z(\tau)})=\chi(\mathfrak{g}_*|_{Z(G)})$ and 
$\chi(\overline{\mathfrak{t}_*}) 
= \chi(\overline{\mathfrak{g}_*})$. 
\end{theorem}

By Theorem \ref{theorem_characterization}, 
the BBK polynomials can be described 
by the train track polynomials and the polynomial 
$\chi(\mathfrak{t}_*|_{\ker \pi_*})$. 
Since these polynomials are conjugacy invariants 
(Theorems \ref{theorem_first-kind} and \ref{theorem_kerpi*_coninv}), 
Theorem \ref{theorem_characterization} 
implies that the BBK polynomials are 
conjugacy invariants of $[f] \in \mathrm{MCG}(\Sigma)$. 

\begin{corollary}
\label{corollary_invariance}
With the same notation as in Theorem \ref{theorem_characterization}, 
the polynomials 
$\chi(\mathfrak{g}_*)$, $\chi(\mathfrak{g}_*|_{Z(G)})$ 
and $\chi(\overline{\mathfrak{g}_*})$ are
conjugacy invariants of $[f] \in \mathrm{MCG}(\Sigma)$.  
\end{corollary}

\begin{remark}
\label{remark_Dodge-Kawamuro}
In 2012, Birman-Brinkmann-Kawamuro \cite{BBK12} stated 
the invariance of the BBK polynomials.  
Dodge-Kawamuro \cite{Kawamuro26} provide more details on the proof 
that the BBK polynomials are invariants. 
They also prove that the BBK polynomials are conjugacy invariants 
of pseudo-Anosov mapping classes. 
However, 
the relation between the train track polynomials and the BBK polynomials is not addressed in \cite{Kawamuro26}. 
\end{remark}

The following result is proved in Section 
\ref{section_conjugacy-invariance} and 
is essentially contained in \cite{HironakaTsang25,BBK12}. 
See Definition \ref{definition_reciprocal} for the notion of 
(anti-)palindromic polynomials.

\begin{proposition}
\label{prop_reciprocal}
Polynomials $\chi(\overline{\mathfrak{t}_*})$, 
$\chi(\mathfrak{t}_*|_{Z(\tau)}) $ and 
$\chi(\mathfrak{t}_*)= \chi(\mathfrak{t}_*|_{Z(\tau)}) 
\chi(\overline{\mathfrak{t}_*})$ are reciprocal, that is, 
palindromic or anti-palindromic. 
\end{proposition}

This paper is primarily based on work done 
as part of a Master's thesis \cite{Mizutani26} by the second author.

The organization of the paper is as follows. 
In Section \ref{section_preliminaries} 
we recall basic definitions and prove some results of train tracks. 
In Section \ref{section_train-track-polynomials} 
we introduce the train track polynomials and 
prove Theorem \ref{theorem_first-kind}. 
In Section \ref{section_BBKpolynomials} 
we introduce the BBK polynomials and review some results of \cite{BBK12}. 
In Section \ref{section_conjugacy-invariance} 
we prove Theorem \ref{theorem_characterization} and  the conjugacy invariance of the BBK polynomials 
in the mapping class groups. 
In Section \ref{section_examples} we give examples.

\section{Preliminaries}
\label{section_preliminaries}

\subsection{Braid groups}
\label{subsection_braid-groups}

Let $B_s$ be the braid group with $s$ strands. 
The group $B_s$ has a presentation 
with generators $\sigma_1, \dots, \sigma_{s-1}$ and relations 
$\sigma_i \sigma_j= \sigma_j \sigma_i$ if $|i-j|>1$ and 
$\sigma_j \sigma_i \sigma_j = \sigma_i \sigma_j \sigma_i$ if $|i-j|=1$.

Let $D_s$ denote a disk with  $s$ punctures $a_1, \dots, a_s$. 
The mapping class group $\mathrm{MCG}(D_s)$ 
is the group of isotopy classes of orientation preserving self-homeomorphisms on $D_s$ preserving  the boundary $\partial D$ 
of the disk setwise.  
There is a surjective homomorphism 
(strictly speaking, anti-homomorphism 
since we read braid words from left to right): 
$$\Gamma: B_s \rightarrow \mathrm{MCG}(D_s)$$
sending the generator $\sigma_i$ to the right-handed half twist $h_i$ 
that interchanges the punctures $a_i$ and $a_{i+1}$ counterclockwise 
for $i= 1, \dots, s-1$.
The kernel of $\Gamma$ is an infinite cyclic group 
generated by a full twist. 
Collapsing the boundary $\partial D$ of $D_s$ into a point  
yields an $s+1$-punctured sphere. 
We may regard $\Gamma(b) \in \mathrm{MCG}(D_s)$ 
as an element of $\mathrm{MCG}(\Sigma_{0, s+1})$.

\subsection{Train tracks}
\label{subsection_tracks}

A {\it train track} $\tau$ is a finite $C^1$ graph 
embedded in a surface $\Sigma$ and 
equipped with a well-defined tangent line at each vertex. 
It also requires that 
each component of the complement $\Sigma \setminus \tau$ is either 
a non-punctured $k$-gon for $k \ge 3$ or 
a once-punctured $k$-gon for $k \ge 1$.  
Edges of a train track are called {\it branches} and vertices  are called 
{\it switches}.  
Two branches of $\tau$ that are tangent at a switch make a {\it cusp}. 
See Figure \ref{fig_local_track}(1).

\begin{definition}
Let $v$ be a switch of $\tau$. 
We label branches 
$b_1', \dots, b_m'$, $b_1, \dots, b_n$ emanating from  $v$ 
that are ordered counterclockwise so that 
the angle between each branch $b_i'$ ($i= 1, \dots, m$) 
and each branch $b_j$ ($j= 1, \dots, n$) is $\pi$ 
(Figure \ref{fig_local_track}(1)(1')). 
We form a {\it tie neighborhood} $\mathcal{N}(\tau) \subset \Sigma$ 
by taking the union of polygons associated with switches of $\tau$ 
and rectangles associated with branches of $\tau$. 
Each rectangle connects the polygons along the branches of $\tau$, 
and both polygons and rectangles are foliated by intervals 
transverse to the branches of $\tau$. 
See  Figure \ref{fig_local_track}(2). 
We call the intervals  {\it ties}. 
\end{definition}

\begin{definition}[Section 3.1 in \cite{HironakaTsang25}]
\label{definition_boundary-component}
Each connected component of $\mathcal{N}(\tau) \setminus \tau$ 
is homeomorphic to an annulus. 
One boundary component of such an annulus is along 
the boundary of $\mathcal{N}(\tau)$ 
and the other boundary component  is along $\tau$. 
The boundary component along $\tau$ is called a 
{\it boundary component of $\tau$}. 
The boundary components of $\tau$ are in one to one correspondence 
with the boundary components of $\mathcal{N}(\tau)$.
We say that a boundary component $c$ of $\tau$ is {\it $k$-pronged} 
if it contains $k$ cusps. 
We say that $c$ is  {\it even-pronged} (resp. {\it odd-pronged})
if it is $k$-pronged for some even (resp. odd) $k$. 
\end{definition}

\begin{definition}
\label{definition_orientable-ttrack}
A train track $\tau$ is {\it orientable} 
if one can orient all branches in such a way that 
at each switch $v$ of $\tau$, the angle between each incoming branch 
toward $v$ 
and each outgoing branch emanating from $v$ is $\pi$.  
If $\tau$ has an odd-pronged boundary component $c$, 
then $\tau$ is non-orientable, 
since there is no consistent orientation on $c$.  
\end{definition}

\subsection{The weight spaces of train tracks} 
\label{subsection_weight-spaces}

\begin{figure}[htbp]
\begin{center}
\includegraphics[width=4.8in]{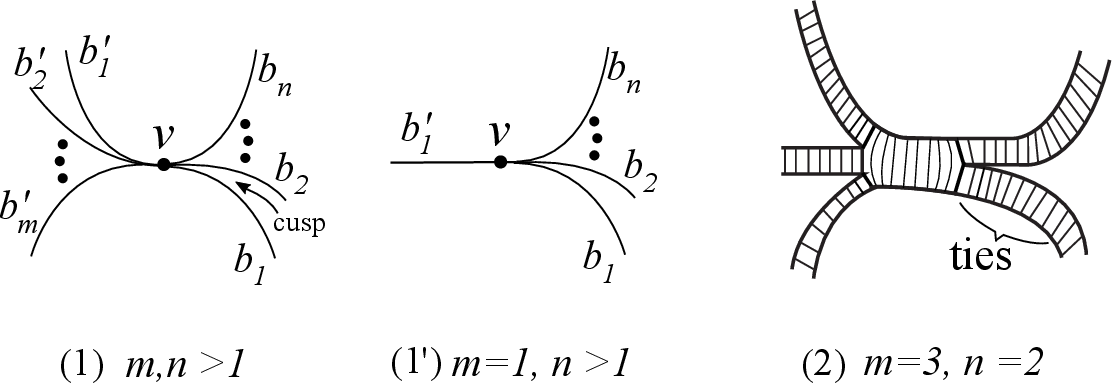}
\caption{(1)(1') Train track $\tau$ near a switch $v$.  
(2) Tie neighborhood $\mathcal{N}(\tau)$.}
\label{fig_local_track}
\end{center}
\end{figure}

 Let $\mathcal{E}(\tau)$ be the set of all branches of $\tau$. 
 Let $V(\tau) \simeq {\Bbb R}^{|\mathcal{E}(\tau)|}$ 
denote the vector space of real weights on the branches of $\tau$. 
We denote  the weight (coefficient) 
that $\eta \in V(\tau)$ assigns 
to  $b \in \mathcal{E}(\tau)$ by $\eta(b)$. 
The {\it weight space} $W(\tau)$ is a subspace of $V(\tau)$ 
consisting of vectors $\eta\in V(\tau)$ 
that satisfy the switch condition 
\begin{center}
$\sum_{i=1}^m \eta(b_i') = \sum_{j=1}^n \eta(b_j)$
\end{center}
at each switch $v$ of $\tau$ (see  Figure \ref{fig_local_track}(1)). 

Let $\sigma$ be a train track in $\Sigma$. 
We say that $\sigma$ is {\it carried by }$\tau$, 
and we write $\sigma< \tau$, 
if $\sigma$ is isotopic to a train track $\sigma'$ 
contained in $\mathcal{N}(\tau)$ and 
$\sigma'$ is transverse to the ties in $\mathcal{N}(\tau)$. 
Let $f: \Sigma \rightarrow \Sigma$ be a homeomorphism. 
A train track $\tau$ is {\it $f$-invariant}  
if $f(\tau) < \tau$. 
In this case, 
collapsing each tie in $\mathcal{N}(\tau)$ into a point,  
we obtain a graph map $\mathfrak{t}: \tau \rightarrow \tau$ 
such  that 
$\mathfrak{t}$ maps vertices to vertices, edges to edge paths. 
Moreover, $\mathfrak{t}$ is locally injective at any points 
that do not map into vertices. 
This map $\mathfrak{t}: \tau \rightarrow \tau$ is called a 
{\it train track representative} of $[f]$.

\begin{definition}[Bestvina-Handel \cite{BestvinaHandel95}]
\label{definition_infinitesimal-real}
A branch $e$ of $\tau$ is said to be {\it infinitesimal} if 
$e$ is eventually periodic under the train track 
representative $\mathfrak{t}: \tau \rightarrow \tau$, that is, 
$\mathfrak{t}^{m+n}(e)= \mathfrak{t}^n(e)$ 
for some integers $m >0$ and $n \ge 0$. 
Other branches of $\tau$ are said to be {\it real}.
\end{definition}

\begin{definition}
\label{definition_transition-matrix}
From a train track representative 
$\mathfrak{t}: \tau \rightarrow \tau$, 
we define a non-negative integer matrix 
$M = (m_{ij})$ 
of degree $|\mathcal{E}(\tau)|$, 
where $(i,j)$-entry $m_{ij}$ is the number of times 
that the image $\mathfrak{t}(e_j)$ of the $j$th edge 
crosses the $i$th edge $e_i$ 
in either direction. 
We call $M$ the  {\it transition matrix} 
of $\mathfrak{t}: \tau \rightarrow \tau$. 
\end{definition}

Suppose that $f: \Sigma \rightarrow \Sigma$ is a pseudo-Anosov map. 
Thurston showed that 
there exists an $f$-invariant train track $\tau$ 
that induces a train track representative 
$\mathfrak{t}: \tau \rightarrow \tau$. 
Moreover, if $M$ is the transition matrix for 
$\mathfrak{t}: \tau \rightarrow \tau$, 
its spectral radius is the stretch factor of $[f]$. 
(See \cite{FLP79} for example.)

The train track representative $\mathfrak{t}: \tau \rightarrow \tau$ 
induces a linear map $\mathfrak{t}_*: V(\tau) \rightarrow V(\tau)$. 
Notice that 
if  $ \eta \in V(\tau)$ that satisfies the switch condition 
at a switch $v$, then 
 $\mathfrak{t}_*(\eta) \in V(\tau)$ also 
satisfies the switch condition at the 
switch $\mathfrak{t}(v)$, since 
$\tau$ is $f$-invariant. 
This implies that 
$\mathfrak{t}_*(W(\tau)) \subset W(\tau)$ 
and we have the linear map 
$\mathfrak{t}_*|_{W(\tau)}: W(\tau) \rightarrow W(\tau)$. 
This restriction map is denoted by 
$$\mathfrak{t}_*: W(\tau) \rightarrow W(\tau).$$
Since the $f$-invariant train track $\tau$ is not unique, 
the characteristic polynomial of 
$\mathfrak{t}_*: V(\tau) \rightarrow V(\tau)$ 
(that is, the characteristic polynomial $\chi(M)$ 
of the above transition matrix $M$) 
is not a well-defined invariant of $[f]$. 
We will verify in Theorem \ref{theorem_inv_t_*} that 
the characteristic polynomial of 
the restriction map 
$\mathfrak{t}_*: W(\tau) \rightarrow W(\tau)$ 
is an invariant of $[f]$.

\subsection{The Thurston symplectic form}
\label{subsection_symplectic-form}

Following \cite{BBK12,HironakaTsang25}, 
we introduce the Thurston symplectic form $\omega$ on $W(\tau)$. 
This is a generalization of the symplectic form on $W(\tau)$ 
for trivalent train tracks \cite{PennerHarer92}. 
Recall that we label branches 
$b_1', \dots, b_m'$, $b_1, \dots, b_n$ emanating from a switch $v$ 
of $\tau$ that are ordered counterclockwise as in 
Figure \ref{fig_local_track}(1)(1'). 
Then the Thurston symplectic form $\omega$ is defined  as follows. 
For $\eta, \zeta \in W(\tau)$, 
  \begin{align*}
  \omega( \eta,\zeta):=\sum_{\substack{\text{switches} \\ \text{in } \tau}}
     \sum_{
      \substack{i < j \\ k < \ell}}
    \left\{
  \begin{vmatrix}
   \eta(b_i) & \eta(b_j) \\
   \zeta(b_i) & \zeta(b_j)
\end{vmatrix}+
\begin{vmatrix}
   \eta(b'_{k}) & \eta(b'_{\ell}) \\
   \zeta(b'_{k}) & \zeta(b'_{\ell})
\end{vmatrix}
\right\},
  \end{align*}
where the second summation is taken over all pairs of branches 
$b_i,b_j$ (resp. $b'_{k}, b'_{\ell}$) 
with $1 \le i < j \le n$ 
(resp. $1 \le k < \ell \le m$). 
In the case $m=1$ 
(Figure \ref{fig_local_track}(1')), 
we ignore 
$\left|\begin{smallmatrix}
\eta(b'_{k}) & \eta(b'_{\ell}) 
\\ \zeta(b'_{k}) & \zeta(b'_{\ell})
\end{smallmatrix}\right|$
and we only take into account of  
$\left|\begin{smallmatrix}
\eta(b_i) & \eta(b_j)
\\ \zeta(b_i) & \zeta(b_j)
\end{smallmatrix}\right|$. 
The case $n=1$ is similar.

\begin{figure}[htbp]
\begin{center}
\includegraphics[width=4.8in]{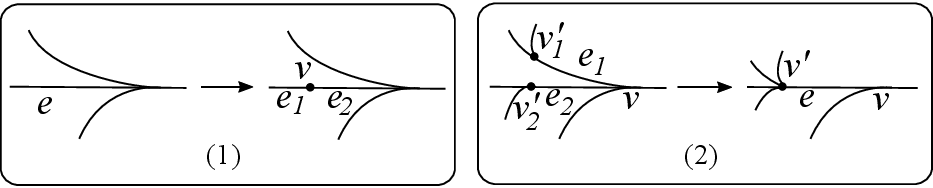}
\caption{(1) Subdivision move. 
(2) Elementary folding move.}
\label{fig_elementary_move}
\end{center}
\end{figure}

We borrow the following definition from \cite{HironakaTsang25}. 
\begin{definition}[Elementary moves] 
Let $\tau \subset \Sigma$ be a train track.

\begin{enumerate}
\item 
Suppose that $e$ is an edge of $\tau$. 
We pick an interior point $v$ of $e$. 
We define a new train track $\tau'$ obtained from $\tau$ 
by replacing $e$ with two new edges $e_1, e_2$ 
whose common endpoint is $v$. 
(Although $\tau$ and $\tau'$ are not isomorphic as a graph, 
they are homeomorphic  as a topological space.)  
We have a natural graph map $\tau \rightarrow \tau'$ 
between the two graphs. 
This map is called the {\it subdivision move}. 
See Figure \ref{fig_elementary_move}(1). 

\item 
Suppose that $e_1$ and $e_2$ are two edges of $\tau$ 
that form a cusp at a switch $v$. 
Let $v_1' \ne v$ and $v_2' \ne v$ be the other endpoints of $e_1$ and 
$e_2$, respectively. 
We define a new train track $\tau'$ obtained from $\tau$ 
by combining $v_1'$ and $v_2'$ into one switch $v'$ 
and replacing $e_1$ and $e_2$ by a single edge $e$ which connects 
$v$ and $v'$. 
There is a graph map $\tau \rightarrow \tau'$ 
sending $v_1$ and $v_2$ to $v$, 
$e_1$ and $e_2$ to $e$ and sending 
the rest of  vertices and edges to themselves. 
This map is called the {\it elementary folding move}. 
See Figure \ref{fig_elementary_move}(2).
\end{enumerate}
We refer to subdivision moves and elementary folding moves 
as {\it elementary moves}.
\end{definition}

\begin{lemma}[Lemmas 4.4 and 4.5 in \cite{HironakaTsang25}]
\label{lem_preserving-omega}
If $m: \tau \rightarrow \tau'$ is an elementary move, then 
$m_*: W(\tau) \rightarrow W(\tau')$ 
is an isomorphism that preserves 
the Thurston symplectic form $\omega$. 
\end{lemma}

The following result follows from a theorem by Stallings \cite{Stallings83}. 

\begin{proposition}
\label{prop_Stalling} 
Let $f: \Sigma \rightarrow \Sigma$ be a pseudo-Anosov map. 
A train track representative 
$\mathfrak{t}: \tau \rightarrow \tau$ of $[f]$  
admits a folding decomposition as follows. 
$$\mathfrak{t}= \rho \circ m_k \circ \cdots \circ m_1$$ 
where 
$m_i: \tau_i \rightarrow \tau_{i+1}$ are  elementary folding moves 
with $\tau_1= \tau_{k+1}= \tau$ and $\rho: \tau \rightarrow \tau$ is a 
homeomorphism.  
\end{proposition}

By Proposition \ref{prop_Stalling} together with 
Lemma \ref{lem_preserving-omega}, 
we  obtain the following. 

\begin{corollary}
\label{cor_preserving-omega}
Let $f: \Sigma \rightarrow \Sigma$ be a pseudo-Anosov map. 
Let $\mathfrak{t}: \tau \rightarrow \tau$ 
be a train track representative of  $[f]$. 
Then $\mathfrak{t}_*: W(\tau) \rightarrow W(\tau)$ 
is an isomorphism that  
preserves the Thurston symplectic form $\omega$. 
\end{corollary}

\subsection{Unstable/stable stars}

The notion of unstable/stable stars  introduced in \cite{HironakaTsang25} plays an important role in the later section.

\begin{figure}[htbp]
\begin{center}
\includegraphics[width=3.8in]{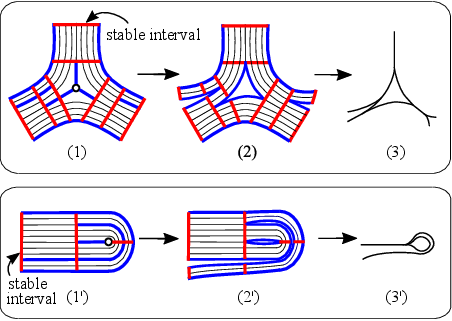}
\caption{
(1)(1') Parts of unstable stars. 
(2)(2') Cutting the surface along unstable stars. 
(3)(3') Train tracks.}
\label{fig_unstable-stars}
\end{center}
\end{figure}

\begin{definition}
Let $f: \Sigma \rightarrow \Sigma$ be a pseudo-Anosov map.
\begin{enumerate}
\item  
Suppose that $f$ is fully punctured. 
Let $x$ be a puncture of $\Sigma$. 
{\it An unstable prong at $x$} is a connected subset 
of an unstable leaf of the unstable foliation $\mathcal{F}^u$ 
that limits to $x$. 
{\it An unstable star $\sigma_x$ at $x$} 
is a maximal disjoint union of unstable prongs at $x$. 
In other words, 
if $\mathcal{F}^u$ is $k$-pronged at $x$, 
an unstable star $\sigma_x$ at $x$ is a disjoint union of $k$ unstable 
prongs at $x$. 
(For example, if $x$ is a puncture of $\Sigma$ 
as shown in Figure \ref{fig_unstable-stars}(1), 
the three unstable prongs at $x$ as in the same figure 
is an unstable star at $x$.)  
The union $\sigma = \bigcup_x \sigma_x$, where $x$ ranges over 
all punctures of $\Sigma$ is called {\it unstable stars of $f$}. 
Stable stars are defined  similarly. 

\item 
Suppose that $f$ is not fully punctured. 
In this case the {\it unstable stars of $f$} 
is the unstable stars of the fully punctured pseudo-Anosov map $f^{\circ}$. 
\end{enumerate}
\end{definition}

Let $\tau \subset \Sigma$ be an  $f$-invariant train track. 
A tie neighborhood $\mathcal{N}(\tau)$ 
can be thought of as the surface 
cut along some unstable stars $ \sigma$ of $f$. 
Conversely, given unstable stars, one can find a train track 
(not necessarily $f$-invariant) with this property:

\begin{lemma}
\label{lemma_ttrack-stars}
Let $f: \Sigma \rightarrow \Sigma$ be a pseudo-Anosov map. 
Let $\sigma$ be unstable stars of $f$. 
Then there is a train track $\tau$  
whose tie neighborhood $\mathcal{N}(\tau)$ can be thought 
of as the surface cut along  $ \sigma \subset \Sigma$. 
\end{lemma}

\begin{proof}
Every stable half-leaf is dense in $\Sigma$. 
(See Lemma 9.6 in \cite{FLP79}.) 
Now we cut $ \Sigma$ along the unstable stars $\sigma$. 
The stable leaves become compact intervals, 
since each stable leaf eventually runs into an unstable star of $f$. 
Then a desired train track can be obtained from the resulting surface 
collapsing each stable interval (corresponding to a tie) 
to a point. 
\end{proof}

We call $\tau$ in  Lemma \ref{lemma_ttrack-stars} 
the {\it train track corresponding to  $\sigma$}. 
The authors learn the following result from Chi Cheuk Tsang.

\begin{proposition}
\label{prop_unstable-stars} 
Let $\sigma$ and $\sigma'$ be unstable stars of a pseudo-Anosov map 
$f: \Sigma \rightarrow \Sigma$. 
Let $\tau$ and $\tau'$ be train tracks 
corresponding to $\sigma$ and $\sigma'$, respectively.  
Then, there exists an elementary folding move 
$\tau \rightarrow \tau'$ 
if and only if $\sigma \supset \sigma'$. 
\end{proposition}

\begin{figure}[htbp]
\begin{center}
\includegraphics[width=3.5in]{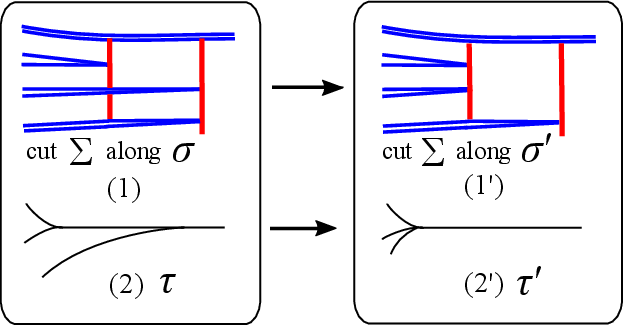}
\caption{Proof of Proposition \ref{prop_unstable-stars}.}
\label{fig_folding-criterion}
\end{center}
\end{figure}

\begin{proof}
(If part.) 
The idea comes from the proof 
\cite[Proposition 3.18]{HironakaTsang25}. 
Suppose that $\sigma \supset \sigma'$. 
By Lemma \ref{lemma_ttrack-stars}, 
tie neighborhoods $\mathcal{N}(\tau)$ and $\mathcal{N}(\tau')$ 
can be thought of as the surfaces cut along unstable stars 
$\sigma$ and $\sigma'$, respectively. 
Since $\sigma$ contains $\sigma'$, 
one can retract $\sigma$ to $ \sigma'$. 
Then the stable intervals (ties) in $\mathcal{N}(\tau)$ 
merge in $\mathcal{N}(\tau')$, 
and this effects an elementary folding move 
$\tau \rightarrow \tau'$.

(Only if part.) 
Recall the construction of train tracks corresponding to 
unstable stars (Figure \ref{fig_unstable-stars}). 
Then it is clear from Figure \ref{fig_folding-criterion} that 
if there is an elementary folding move 
$\tau \rightarrow \tau'$, then 
$\sigma$ contains $\sigma'$. 
\end{proof}

\begin{corollary}
\label{cor_f-invariant}
Let $f: \Sigma \rightarrow \Sigma$ be a pseudo-Anosov map. 
Let $\sigma$ be unstable stars of $f$. 
Suppose that $\tau$ is the train track 
corresponding to $\sigma$. 
Then 
$f(\tau)< \tau$ (i.e., $f$-invariant) 
if and only if $f(\sigma) \supset \sigma$. 
\end{corollary}

\begin{proof}
The statement follows from Propositions \ref{prop_Stalling} 
and \ref{prop_unstable-stars}. 
\end{proof}

\subsection{The radical of the Thurston symplectic form}

The Thurston symplectic form $\omega$ is not a symplectic form 
in general, that is, 
the {\it radical} $Z(\tau)$ of $\omega$ 
 is not  trivial in general.  
 $$Z(\tau)=\{ \eta \in W(\tau) \mid 
  \omega(\eta,\zeta) =0 
 \text{ for all }  \zeta \in W(\tau)  \}.$$

\begin{figure}[htbp]
\begin{center}
\includegraphics[width=4.8in]{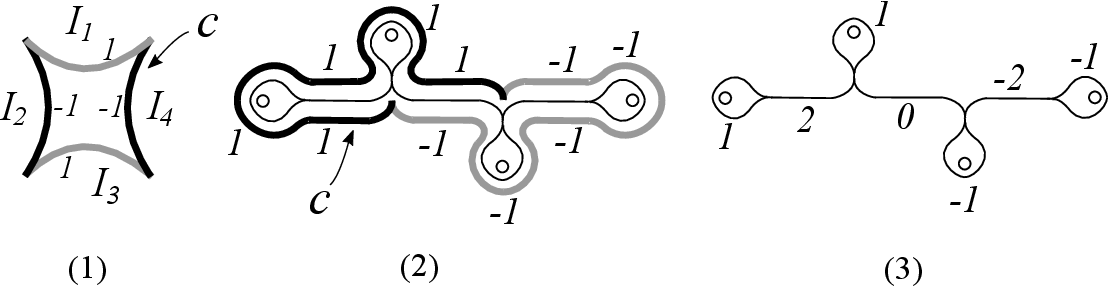}
\caption{(1)$4$-pronged boundary component.
(2) $2$-pronged boundary component of $\tau$.
(3) Radical element $r_c$.}
\label{fig_radical-element}
\end{center}
\end{figure}

\begin{definition}[The radical element  \cite{HironakaTsang25}]
\label{definition_rc} 
Let $c$ be an even-pronged boundary component of $\tau$.
We label the intervals on $c$ between consecutive cusps by
$I_1, \dots, I_{2n}$ in  cyclic order. 
(Hence, the endpoints of each $I_i$ are cusps of $\tau$.)
The {\it radical element} $r_c$  is an element of $W(\tau)$
defined as follows. 
First, assign weight $(-1)^i$ to each edge of $\tau$
that lies on $I_i$ (Figure \ref{fig_radical-element}(1)(2)), 
and assign weight $0$ 
to each edge disjoint from all the intervals $I_i$.
Then, for each edge $e$ of $\tau$, 
define the coefficient $r_c(e)$ 
to be the sum of weights assigned to $e$ 
over all its occurrences in the intervals 
(Figure \ref{fig_radical-element}(2)(3)).
In other words, if $e$ appears multiple times, possibly in multiple intervals $I_i$,
then $r_c(e)$ is the sum of  weights 
assigned to $e$ at each occurrence.
In particular, if an edge $e$ appears exactly once in exactly one of the intervals $I_i$, then $r_c(e)=(-1)^i$.
It is proved in \cite[Proposition 5.2]{HironakaTsang25} that 
 $r_c$ is an element of the radical $Z(\tau) \subset W(\tau)$. 
\end{definition}

Suppose that $\tau'$ is obtained from $\tau$ by 
an elementary folding move. 
Then there is a natural one to one correspondence 
between the boundary components of $\tau$ and those of $\tau'$. 
The following result can be found in 
 \cite[Proposition 5.3]{HironakaTsang25}. 

\begin{proposition}
\label{prop_prop5.3}
Let $m: \tau \rightarrow \tau'$ be an elementary move. 
Suppose that $c'$ is a boundary component of $\tau'$ 
that corresponds to an even-pronged boundary component $c$ of $\tau$. 
Then the isomorphism 
$m_*: W(\tau) \rightarrow W(\tau')$ sends 
$r_c$ to $r_{c'}$. 
\end{proposition}

\begin{definition}
Let $l$ be a lamination in $\Sigma$. 
Let $\tau \subset \Sigma$ be a train track. 
We say that $\tau$ {\it carries} $l$ 
if leaves of $l$ are contained in a  tie neighborhood $\mathcal{N}(\tau)$ 
and transverse to the ties in $\mathcal{N}(\tau)$.  
We say that $\tau$ {\it fully carries} $l$ 
if $\tau$ carries $l$ and every tie in $\mathcal{N}(\tau)$ 
intersects some leaf of $l$. 
\end{definition}

If $\tau$ is a train track corresponding to some unstable stars of 
a pseudo-Anosov map, then 
$\tau$ fully carries the unstable lamination of $f$. 

The following definition is used in the proof of Proposition \ref{prop_prop8.1}. 

\begin{definition}[Definition 4.7 in \cite{HironakaTsang25}]
\label{definition_equivalence-track}
We consider the equivalence class on the set of all train tracks 
that is generated by there being an elementary move between 
two train tracks. 
We say that train tracks $\tau$ and $\tau'$ are {\it equivalent} 
if $\tau$ and $\tau'$ are in the same equivalence class. 
\end{definition}

The following result is a refinement of 
Proposition 8.1 in \cite{HironakaTsang25}
and the proof below is due to Chi Cheuk Tsang.

\begin{proposition}
\label{prop_prop8.1}
Let $f: \Sigma \rightarrow \Sigma$ be a pseudo-Anosov map, and 
let $\sigma$ be  unstable stars of $f$. 
Then the train track $\tau$ corresponding to $\sigma$ 
satisfies the following. 
\begin{enumerate}
\item 
The radical $Z(\tau)$ is given by 
$$Z(\tau) = \mathrm{span}\{r_c\},$$ 
where $c$ ranges over all even-pronged boundary components of $\tau$. 

\item 
The radical elements $\{r_c\}$ 
are linearly independent if and only if $\tau$ is non-orientable. 

\item 
The dimension of $Z(\tau)$ is 
$$
\begin{array}{ll}
\# \{\mbox{even-pronged boundary components}\}-1
& \mbox{if}\ \tau \  \mbox{is orientable},  
\\
\# \{\mbox{even-pronged boundary components}\} 
& \mbox{if}\ \tau\ \mbox{is non-orientable}. 
\end{array}
$$
\end{enumerate}
\end{proposition}

We will need  {\it standardly embedded train tracks}  
introduced by Hironaka-Tsang. 
See Definition 3.1 of \cite{HironakaTsang25}. 
The property we need on standardly embedded train tracks $\tau$
is that 
one side of a switch of $\tau$ consists of a number of real edges 
(possibly one real edge) and the other side consists of 
two infinitesimal edges on infinitesimal polygons.  
See Figure \ref{fig_standardly}(1)(2). 
For example, 
the train track in Figure \ref{fig_radical-element}(3) is 
a standardly embedded train track.

\begin{proof}[Proof of Proposition \ref{prop_prop8.1}]
The claim (1) holds 
from Proposition 8.1 in \cite{HironakaTsang25}. 

We turn to (2). 
For the proof of `the only if part' ($\Longrightarrow$), 
we prove the contrapositive. 
Suppose that $\tau$ is orientable. 
Then we can consistently orient all edges of $\tau$, and 
all boundary components of $\tau$ are even-pronged. 
(See Definition \ref{definition_orientable-ttrack}.) 
Let $c$ be any boundary component of $\tau$ and 
let $I_1, \dots, I_{2n}$ be the intervals on $c$ 
as in Definition \ref{definition_rc}. 
We  fix a sign of the radical element $r_c$ 
(i.e. signs of coefficients of edges in $r_c$) 
by specifying 
that $I_{2i+1}$ are oriented clockwise and 
$I_{2i}$ are oriented counterclockwise. 
Then $ \sum_c r_c=0$, 
since each edge of $\tau$ appears twice in the sum, 
once on a clockwise interval and the other time on a counterclockwise 
interval. 
(If $\tau$ has a single boundary component $c$, 
then $r_c=0$.) 
Hence $\{r_c\}$ are linearly dependent.

We turn to `the if part' ($\Longleftarrow$) of (2). 
Suppose that $\tau$ is a non-orientable train track. 
We prove that $\{r_c\}$ are linearly independent. 
To do this, we first consider the case 
where $\tau$ has at least two boundary components. 
The validity of the claim is invariant under equivalence of 
train tracks (see Definition \ref{definition_equivalence-track}). 
We note that 
any two train tracks $\tau$ and $\tau'$ corresponding to unstable 
stars $ \sigma$ and $\sigma'$ are equivalent. 
Hence, we can take particular unstable stars $\sigma$ 
and the corresponding train rack $\tau$. 
From the proof of Proposition 5.4 in \cite{HironakaTsang25}, 
 we can assume that 
all boundary components of $\tau$, 
except a specified boundary component $O$, 
are infinitesimal polygons, 
that is, the polygons whose sides are infinitesimal edges. 
(In fact, choose unstable stars $\sigma$ so that 
the length of the unstable star at a puncture corresponding to $O$ 
is long enough. Then the corresponding train track $\tau$ has 
the above property.) 
We can further assume that 
$\tau$ is a standardly  embedded train track 
(c.f. Construction 3.6 of \cite{HironakaTsang25}). 
If $\tau$ has an odd-pronged boundary component, 
then we can choose $O$ to be odd-pronged. 
In this case, all even-pronged boundary components are disjoint. 
This means that $\{r_c\}$ are linearly independent. 

We are left with the case where 
$\tau$ is non-orientable and 
all boundary components are even-pronged. 
We assume that there is a relation $\sum_c a_c r_c=0$. 
Since the boundary components other than $O$ are disjoint, 
$\{r_ c\}_{c \ne O}$, 
that is, 
the radical elements other than $r_O$, 
are linearly independent, and hence 
the coefficient $a_O$ is nonzero. 
We claim that 
one can construct orientations of edges on $ \tau$ 
using the relation $\sum_c a_c r_c=0$ 
and this would contradicts the assumption 
that $\tau$ is non-orientable.

To prove this claim, 
we choose the labeling $I_1, \dots, I_ {2n}$ 
on  $O$ as in Definition \ref{definition_rc}. 
Then orientations on the edges on $O$ are defined by 
orientating all edges on $I_{2i+1}$ 
clockwise and all edges on $I_{2i}$ counterclockwise. 
We first verify that  this is well defined: 
The only way this can fail is that 
an edge appears twice on $O$, on $I_{2i}$ and $I_{2j}$, 
or $I_{2i+1}$ and $I_{2j+1}$. 
If this is the case, 
the coefficient of the edge in the vector  $\sum_c a_c r_c$ 
must be $\pm 2a_O$ 
which is nonzero. 
This contradicts the relation $\sum_c a_c r_c=0$.

We continue with the proof of the claim. 
We check that the orientations induced by $O$ matches up 
at each switch of $\tau$: 
Since $\tau$ is a standardly embedded train track, 
one side of a switch of $\tau$ consists of a number of real edges 
and the other side consists of  two infinitesimal edges, 
say $e_1$ and $e_2$. 
In the former case, 
each pair of adjacent real edges is consistently oriented 
as we proved above. 
In the latter case, 
we suppose that 
the two infinitesimal edges $e_1$ and $e_2$ 
belong to an infinitesimal polygon $C$. 
We look at the relation $\sum_c a_c r_c=0$. 
Then $e_1$ and $e_2$ meet the two boundary components of 
$\tau$, namely $O$ and $C$. 
Hence, we have the two equations coming from the coefficients of 
$e_1$ and $e_2$  of the form 
$\pm a_O \pm a_C=0$, 
since the coefficients of $e_1$ and $e_2$ in 
the two radical elements $r_O$ and $r_C$ are  $1$ or $-1$. 
Note that the coefficients of $e_1$ and $e_2$ 
in $r_C$ have the opposite signs, 
since they are adjacent on $C$. 
Hence, the coefficients of $e_1$ and $e_2$ in $r_O$ 
have the opposite signs, 
for otherwise $a_O= a_C= 0$. 
This implies that 
one infinitesimal edge, say $e_1$ lies on $I_{2i}$ (with even index)
of the side on $O$,  
while the other infinitesimal edge $e_2$  lies on $I_{2j+1}$ 
(with odd index) of the side on $O$ 
(see Figure \ref{fig_standardly}(3)). 
Therefore, the orientations induced by $O$ 
on $e_1$ and $e_2$ match up. 
This tells us that $\tau$ is orientable, which is a contradiction.

\begin{figure}[htbp]
\begin{center}
\includegraphics[width=4.4in]{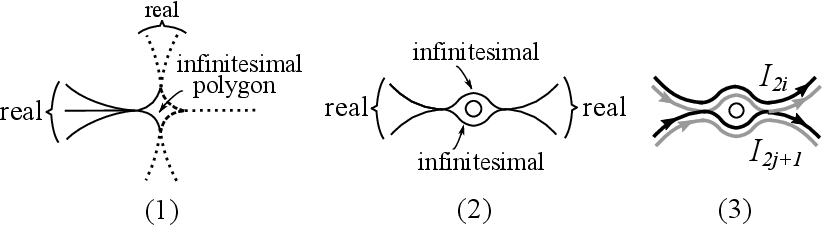}
\caption{(1)(2) Standardly embedded train tracks near switches. 
(3) Proof of Proposition \ref{prop_prop8.1}(2).   
(The symbol $\circ$ denotes a puncture of the surface.)}
\label{fig_standardly}
\end{center}
\end{figure}

We next consider the remaining case 
where $\tau$ is a non-orientable train track 
with a single even-pronged boundary component $O$. 
Choose the labeling $I_1, \dots, I_ {2n}$ 
on  $O$ as in Definition \ref{definition_rc}. 
The orientations on the edges of $O$ (hence on the edges of $\tau$) 
are defined by orientating all edges on $I_{2i+1}$ 
clockwise and all edges on $I_{2i}$ counterclockwise. 
By a similar argument as above, 
the relation $r_O=0$ implies that 
one can consistently orient all edges of $\tau$, 
which means $\tau$ is orientable. 
This contradicts the assumption that $\tau$ is non-orientable. 
This completes the proof of (2). 

The claim (3) follows from (1) and (2) together with 
the proof of (2). 
\end{proof}

\section{Train track polynomials}
\label{section_train-track-polynomials}

In this section, we introduce the first kind of polynomials, 
which we call the train track polynomials, 
and prove that 
they are conjugacy invariants of pseudo-Anosov mapping classes. 
We will need the following lemma.

\begin{lemma}
\label{lemma_Rykken}
Let $A,B$ and $C$ be vector spaces. 
Suppose that we have the commutative diagram 
\[
\begin{tikzcd}
0 \arrow[r]
  & A \arrow[r,"f"] \arrow[d,"\alpha"]
  & B \arrow[r,"g"] \arrow[d,"\beta"]
  & C \arrow[r] \arrow[d,"\gamma"]
  & 0 \\
0 \arrow[r]
  & A \arrow[r,"f"]
  & B \arrow[r,"g"]
  & C \arrow[r]
  & 0
\end{tikzcd}
\]
where the rows are exact and 
$\alpha, \beta$ and $\gamma$ are representation matrices 
corresponding to the linear transformations. 
Then there exists a choice of basis such that 
$\beta= 
\left(\begin{array}{cc}\alpha & * \\ 
0 & \gamma\end{array}\right)$. 
In particular, 
the characteristic polynomial of $\beta$ is factorized 
as 
$\chi(\beta)= \chi(\alpha) \chi(\gamma)$. 
\end{lemma}
The proof of the lemma is straightforward.

\begin{theorem}
\label{theorem_inv_t_*} 
Let $f: \Sigma \rightarrow \Sigma$ be a pseudo-Anosov map. 
Let $\mathfrak{t}:  \tau \rightarrow \tau$ be a train track representative of $[f]$. 
Then the characteristic polynomial $\chi(\mathfrak{t}_*)$ of 
$\mathfrak{t}_*: W(\tau) \rightarrow W(\tau)$ 
is a conjugacy invariant of $[f] \in \mathrm{MCG}(\Sigma)$. 
\end{theorem}

\begin{proof}
Suppose that 
$\mathfrak{t}: \tau \rightarrow \tau$ and 
$\mathfrak{t}': \tau' \rightarrow \tau'$ 
are train track representatives of $[f]$. 
Let $\mathfrak{t}_*: W(\tau) \rightarrow W(\tau)$ and 
$\mathfrak{t}'_*: W(\tau') \rightarrow W(\tau')$ 
be the linear maps. 
Tie neighborhoods of $\tau$ and $\tau'$ 
can be thought of as the surfaces  
cut along unstable stars, say $\sigma$ and $\sigma'$, respectively 
of the common unstable foliation of $f$. 
Since both $\tau$ and $\tau'$ are $f$-invariant,  
we have 
$f(\sigma) \supset \sigma$ and $f(\sigma') \supset \sigma'$ 
by Corollary \ref{cor_f-invariant}. 
We set $\sigma'' := \sigma \cap \sigma'$. 
Then the unstable stars $\sigma'' $ also satisfy 
$f(\sigma'') \supset \sigma''$.
This tells us that 
the train track $\tau''$ corresponding to  $\sigma''$ is also 
$f$-invariant 
by Corollary \ref{cor_f-invariant}. 
Hence, we have a train track representative 
$\mathfrak{t}'': \tau'' \rightarrow \tau''$ of $[f]$.

By definition of the unstable stars $\sigma''= \sigma \cap \sigma'$, 
there exist homotopy equivalence graph maps 
$m: \tau \rightarrow \tau''$ and $m': \tau' \rightarrow \tau''$. 
By a result of Stallings \cite{Stallings83}, each of the graph maps 
$m: \tau \rightarrow \tau''$ and $m': \tau' \rightarrow \tau''$ 
can be factored as a composition of 
elementary folding moves. 
Then the following commutative diagram holds. 
\begin{equation}
\label{equation_commutative-1}
\begin{tikzcd}
    W(\tau) \arrow[r, "m_*"] 
    \arrow[d, "\mathfrak{t}_* "]&W(\tau'') 
    \arrow[d,"\mathfrak{t}''_* "]\\
    W(\tau) \arrow[r, "m_*"] & W(\tau'')
 \end{tikzcd}
\end{equation}
where 
$m_*:  W(\tau) \rightarrow W(\tau'')$ 
is the isomorphism induced by 
$m$ (see Lemma \ref{lem_preserving-omega}). 
To see  the commutative square  
$m_* \circ \mathfrak{t}_*= \mathfrak{t}'' \circ m_*$, 
we note that the train track representatives $\mathfrak{t}$ and $\mathfrak{t}''$ 
can be factorized into elementary folding moves 
(see Proposition \ref{prop_Stalling}).   
For each  element $\eta \in W(\tau)$, we keep track of the weights 
on the branches, and add them up 
when we perform elementary folding moves. 
Then we obtain the desired claim   
$m_* \circ \mathfrak{t}_* (\eta)= \mathfrak{t}'' \circ m_* (\eta)$ 
for $\eta \in W(\tau)$. 

By (\ref{equation_commutative-1}), we have 
$\chi(\mathfrak{t}_*)= \chi(\mathfrak{t}''_*)$. 
Applying the other  homotopy equivalence graph map 
$m': \tau' \rightarrow \tau''$ to the same argument, 
we have $\chi(\mathfrak{t}'_*)= \chi(\mathfrak{t}''_*)$. 
Hence $\chi(\mathfrak{t}_*) = \chi(\mathfrak{t}'_*)$ 
is an invariant of $[f]$. 

We turn to the conjugacy invariance of $\chi(\mathfrak{t}_*)$. 
Suppose that pseudo-Anosov maps 
$f,f': \Sigma \rightarrow \Sigma$ are conjugate, 
that is, there exists an orientation preserving homeomorphism 
$u: \Sigma \rightarrow \Sigma$ such that 
$f'= u \circ f \circ u^{-1}$. 
Let $\mathfrak{t}: \tau \rightarrow \tau$ 
be a train track representative of $[f]$. 
We set $\tau':= u(\tau)$. 
Then the train track $f'(\tau')=  u \circ f \circ u^{-1}(\tau')$ 
is carried by $\tau'$, and 
we obtain a train track representative 
$\mathfrak{t}': \tau' \rightarrow \tau'$ of 
$[f'] \in \mathrm{MCG}(\Sigma)$. 
Then the following commutative diagram holds. 
$$
\begin{tikzcd}
    W(\tau) \arrow[r, "u_* "] 
    \arrow[d, "\mathfrak{t}_*"]&W(\tau') 
    \arrow[d,"\mathfrak{t}_*'"]\\
    W(\tau) \arrow[r, "u_* "] & W(\tau')
 \end{tikzcd}
$$
where $u_*: W(\tau) \rightarrow W(\tau')$ 
is the isomorphism induced by $u|_{\tau}: \tau \rightarrow \tau'$. 
Hence, 
$\chi(\mathfrak{t}_*) = \chi(\mathfrak{t}'_*)$ 
is a conjugacy invariant of $[f]$. 
\end{proof}

By Corollary \ref{cor_preserving-omega} 
we have the restriction map 
$$\mathfrak{t}_*|_{Z(\tau)}: Z(\tau) \rightarrow Z(\tau).$$ 
The quotient space 
$W(\tau)/Z(\tau)$ inherits the Thurston symplectic form $\omega$ 
and is a non-degenerate symplectic vector space.  
The induced map $\overline{\mathfrak{t}_*}$ on $W(\tau)/Z(\tau)$ 
preserves the form $\omega$. 
Then, it is a well-known fact that 
$\chi (\overline{\mathfrak{t}_*})$ is palindromic or anti-palindromic. 
The following commutative diagram holds. 
$$
\begin{tikzcd}
0 \arrow[r]
  & Z(\tau) \arrow[r] \arrow[d,"\mathfrak{t}_*|_{Z(\tau)}"]
  & W(\tau) \arrow[r] \arrow[d,"\mathfrak{t}_*"]
  & W(\tau)/Z(\tau) \arrow[r] \arrow[d,"\overline{\mathfrak{t}_*}"]
  & 0 \\
0 \arrow[r]
  & Z(\tau) \arrow[r]
  & W(\tau) \arrow[r]
  & W(\tau)/Z(\tau) \arrow[r]
  & 0
\end{tikzcd}
$$
By the above diagram together with 
Lemma \ref{lemma_Rykken}, one has the following. 

\begin{proposition}
\label{prop_tt-homology-polynomial}
The characteristic polynomial of 
$\mathfrak{t}_*: W(\tau) \rightarrow W(\tau)$ has a factorization  
$\chi(\mathfrak{t}_*)= \chi(\mathfrak{t}_*|_{Z(\tau)}) 
\chi(\overline{\mathfrak{t}_*})$. 
\end{proposition}

We call $\chi(\mathfrak{t}_*|_{Z(\tau)})$ and 
$\chi(\overline{\mathfrak{t}_*})$ 
the {\it radical polynomial} and 
the {\it (train track) symplectic polynomial} of $[f]$, respectively. 
We collectively refer to 
$\chi(\mathfrak{t}_*|_{Z(\tau)})$, 
$\chi(\overline{\mathfrak{t}_*})$ 
and 
$\chi(\mathfrak{t}_*)= \chi(\mathfrak{t}_*|_{Z(\tau)}) 
\chi(\overline{\mathfrak{t}_*})$ 
as the {\it train track polynomials}.

\begin{theorem}
\label{theorem_invariant-tt-puncture}
The radial polynomial $\chi(\mathfrak{t}_*|_{Z(\tau)})$ 
and the train track symplectic polynomial
$\chi(\overline{\mathfrak{t}_*})$ 
are 
conjugacy invariants of $[f] \in \mathrm{MCG}(\Sigma)$. 
\end{theorem}

\begin{proof}
By Proposition \ref{prop_tt-homology-polynomial}, 
$\chi(\overline{\mathfrak{t}_*})$ is given by 
$\chi(\overline{\mathfrak{t}_*})= \frac{\chi(\mathfrak{t}_*)}{\chi(\mathfrak{t}_*|_{Z(\tau)})}$. 
By Theorem \ref{theorem_inv_t_*}, 
it is enough to prove that 
$\chi(\mathfrak{t}_*|_{Z(\tau)})$ 
is a conjugacy invariant of $[f]$.

Suppose that 
$\mathfrak{t}: \tau \rightarrow \tau$ and 
$\mathfrak{t}': \tau' \rightarrow \tau'$ 
are train track representatives of $[f]$. 
We take a third train track $\tau''$ 
and a train track representative 
$\mathfrak{t}'': \tau'' \rightarrow \tau''$ 
as in the proof of Theorem \ref{theorem_inv_t_*}. 
Recall that the homotopy equivalence graph map 
$m: \tau \rightarrow \tau''$ is factored as a composition of 
elementary folding moves. 
Lemma \ref{lem_preserving-omega} 
tells us that 
the isomorphism 
$m_*: W(\tau) \rightarrow W(\tau'')$ 
preserves the Thurston symplectic form $\omega$, 
and we have the restriction map 
$m_*|_{Z(\tau)}: Z(\tau) \rightarrow Z(\tau'')$. 
This together with 
(\ref{equation_commutative-1}) yields 
the following commutative diagram. 
\begin{equation}
\label{equation_com_Z}
\begin{tikzcd}
    Z(\tau) \arrow[r, "m_*|_ {Z(\tau)}"] 
    \arrow[d, " \mathfrak{t}_*|_{Z(\tau)}"]
    &Z(\tau'') \arrow[d,"\mathfrak{t}''_*|_{Z(\tau'')}"]\\
    Z(\tau) \arrow[r, " m_*|_{Z(\tau)} "] & Z(\tau'')
 \end{tikzcd}
\end{equation}
Hence, we have 
$\chi(\mathfrak{t}_*|_{Z(\tau)}) = \chi(\mathfrak{t}''_*|_{Z(\tau'')})$. 
Applying the same argument to the other 
homotopy equivalence graph map $m': \tau' \rightarrow \tau''$, 
we obtain 
$\chi(\mathfrak{t}'_*|_{Z(\tau')})=\chi(\mathfrak{t}''_*|_{Z(\tau'')})$. 
Thus, $\chi(\mathfrak{t}_*|_{Z(\tau)})=\chi(\mathfrak{t}'_*|_{Z(\tau')})$ 
is an invariant of $[f]$. 
One can prove the conjugacy invariance of 
$\chi(\mathfrak{t}_*|_{Z(\tau)})$ 
in the same manner as in the proof of Theorem \ref{theorem_inv_t_*}. 
\end{proof}

\begin{proof}[Proof of Theorem \ref{theorem_first-kind}]
See Theorems \ref{theorem_inv_t_*} and 
\ref{theorem_invariant-tt-puncture}. 
\end{proof}

The train track polynomials are defined directly 
from  train track representatives. 
Pseudo-Anosov maps  on different surfaces 
may have the same train track representative: 
Let $f: \Sigma \rightarrow \Sigma$ be a pseudo-Anosov map. 
A singularity $q$ of the unstable foliation $\mathcal{F}^u$ 
is either unpunctured or punctured. 
Suppose that $q$ is an unpunctured  singularity. 
Removing the periodic orbit $\mathcal{O}(q)$  under the action of $f$, 
we obtain a surface $\Sigma':= \Sigma \setminus \mathcal{O}(q)$ 
and a new pseudo-Anosov map 
$f':= f|_{\Sigma'}: \Sigma' \rightarrow \Sigma'$. 
Then 
a train track representative $\mathfrak{t}: \tau \rightarrow \tau$ 
of  $[f] \in \mathrm{MCG}(\Sigma)$ 
is also a train track representative of 
$[f'] \in \mathrm{MCG}(\Sigma')$. 
This implies that the train track polynomials of $[f]$ and $[f']$ 
are the same:  

\begin{proposition}
\label{prop_same-ttp} 
The radical polynomial and the (train track) symplectic polynomial 
of $[f]$ are the same as  those of $[f']$. 
\end{proposition}

The train track polynomials cannot distinguish 
between the two mapping classes $[f]$ and $[f']$. 
In Section \ref{section_BBKpolynomials}, 
we introduce the second kind of polynomials, the BBK polynomials. 
We will prove in Section \ref{section_conjugacy-invariance} that 
they have a finer structure than the train track polynomials.

\section{BBK polynomials} 
\label{section_BBKpolynomials}

\subsection{Bestvina-Handel algorithm}

Let $G$ be a finite graph,  
and $\mathcal{V}(G)$ the set of vertices of $G$. 
We denote by $\mathcal{E}_{\mathrm{dir}}(G)$ and $\mathcal{E}(G)$, 
the set of oriented  edges of $G$ and 
the set of unoriented edges of $G$, respectively.
If $e \in \mathcal{E}_{\mathrm{dir}}(G)$, 
we use the notation $\overline{e}$ to denote the edge $e$ 
with the reversed orientation.  
The {\it link} $\mathrm{Lk}(v) \subset \mathcal{E}_{\mathrm{dir}}(G)$ 
is the set of oriented edges emanating from a vertex $v$ of $G$.

\begin{figure}[htbp]
\begin{center}
\includegraphics[width=3.5in]{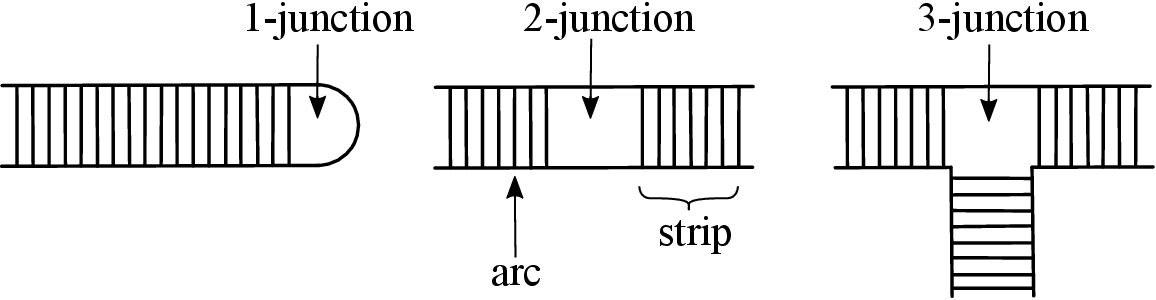}
\caption{Fibered surface.}
\label{fig_fiberedsurface}
\end{center}
\end{figure}

A {\it fibered surface} $F$ is a compact surface 
with boundary which is decomposed into arcs and polygons 
modeled on $k$-junctions 
where $k= 1,2,\dots$ (Figure \ref{fig_fiberedsurface}).  
An arc or a junction is called a {\it decomposition element}. 
Components fibered by arcs are called {\it strips}. 
Collapsing each arc to a point and each junction to a point 
gives a graph $G$ whose vertices come from junctions 
and edges come from strips. 
We can regard $G$ as a graph embedded in $F$ 
and $G$ represents the spine of $F$.

Let $\Sigma= \Sigma_{g,s}$ be a punctured surface for some $s \ge 1$. 
Consider a homeomorphism $f: \Sigma \rightarrow \Sigma$. 
We say that a fibered surface $F \subset \Sigma$ 
{\it carries} $f$ if the inclusion $F \rightarrow \Sigma$ 
is a homotopy equivalence, 
$f$ maps each decomposition element of $F$ 
into a decomposition element of $F$, 
and $f$ maps each junction into a junction. 
(In particular, $f$ maps  $F$ into itself.) 
If $F $ carries $f$, 
collapsing each decomposition element of $F$ to a point, 
$f$ determines a graph map  
$\mathfrak{g}: G \rightarrow G$ 
on the spine $G$ of $F$ 
sending vertices to vertices, and 
each edge to an edge-path. 
We call $\mathfrak{g}: G \rightarrow G$ 
an {\it induced graph map} of $[f] \in \mathrm{MCG}(\Sigma)$. 
From such a graph map $\mathfrak{g}: G \rightarrow G$, 
we define the transition matrix $T$ of degree $|\mathcal{E}(G)|$ 
in the same manner as in Definition \ref{definition_transition-matrix}. 
The spectral radius of $T$ is called 
the {\it growth rate} of $\mathfrak{g}$. 
The graph map $\mathfrak{g}: G \rightarrow G$ 
is {\it efficient} if for each $k \ge 1$, 
the graph map $\mathfrak{g}^k: G \rightarrow G$ 
has no back track, i.e., 
$\mathfrak{g}^k$ maps each oriented edge of $G$ to an edge path 
that does not contain an edge $e \in \mathcal{E}_{\mathrm{dir}}(G)$ 
followed by $\overline{e}$.

Since $\Sigma$ is a punctured surface, 
we can always find a fibered surface that carries $f$ 
and we have an induced graph map 
$\mathfrak{g}: G \rightarrow G$ 
(which is not necessarily efficient). 
The Bestvina-Handel algorithm proceeds as follows. 
If the initial induced graph map $\mathfrak{g}: G \rightarrow G$  
is not efficient, then 
either there is a $\mathfrak{g}$-invariant proper subgraph of $G$ 
(corresponding to reducing curves for the reducible mapping class $[f]$), 
or it is possible to replace $G$ with a graph $G'$ such that 
the inclusion $G' \rightarrow \Sigma$ is a homotopy equivalence and 
the induced graph map $\mathfrak{g}': G' \rightarrow G'$ 
has a smaller growth rate than $\mathfrak{g}: G \rightarrow G$. 
There is a number $N$ depending on the initial graph map 
$\mathfrak{g}: G \rightarrow G$ and the surface $\Sigma$ 
such that, after at most $N$ steps, 
either a reduction or an efficient graph map 
$\mathfrak{g}: G \rightarrow G$ can be found, 
where we use the same notation $\mathfrak{g}: G \rightarrow G$ 
for the resulting graph map as  the initial one. 
Such an efficient graph map can be thought of 
as a simplest possible 
induced graph map of $[f]$ in the sense of the growth rate. 
Once an efficient graph map $\mathfrak{g}: G \rightarrow G$ is found, 
one can determine the Nielsen-Thurston type of $f$. 
In the pseudo-Anosov case, 
we call the resulting graph $G$ a {\it Bestvina-Handel graph} 
and $\mathfrak{g}: G \rightarrow G$ a {\it train track map} of $[f]$. 
In this case, the growth rate of $\mathfrak{g}$ is equal to 
the stretch factor of  $[f]$. 

In the case where $\Sigma$ is an $s$-punctured disk $D_s$, 
the program ``Trains'' 
developed by Toby Hall is available \cite{Hall}. 
Using this program. one can easily find train track maps 
of pseudo-Anosov braids.

\begin{figure}[htbp]
\begin{center}
\includegraphics[width=4in]{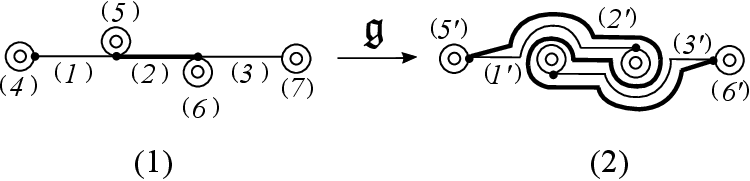}
\caption{Example \ref{example_4braid_BH}. 
The label (i') means the image of the edge labeled (i).}
\label{fig_g_1-1-1}
\end{center}
\end{figure}

\begin{example}
\label{example_4braid_BH}
Let $\Gamma: B_s \rightarrow \mathrm{MCG}(D_s)$ be 
the surjective homomorphism as in Section \ref{subsection_braid-groups}. 
Consider 
$[f]= \Gamma(\sigma_1 \sigma_3 \sigma_2^{-1}) \in \mathrm{MCG}(D_4)$. 
This is a pseudo-Anosov mapping class. 
The graph map $\mathfrak{g}: G \rightarrow G$ in Figure~\ref{fig_g_1-1-1} 
is  a train track map of $[f]$, 
where the edges are labeled $(1),\dots, (7)$.  
The transition matrix of $\mathfrak{g}$ is given by 
$T= 
\left(\begin{array}{cc}M & 0 \\ 
* & Q\end{array}\right)$ 
where $M = 
\left(\begin{smallmatrix}
1 & 1 & 0 
\\ 1 & 3 & 1 
\\ 0 & 1 & 1
\end{smallmatrix}\right)$ 
and 
$Q=\left(\begin{smallmatrix}
0 & 1 & 0 & 0
\\ 0 & 0 & 0 & 1 
\\ 1 & 0 & 0 & 0
\\ 0 & 0& 1& 0
\end{smallmatrix}\right)$. 
Then $\chi(T)= \chi(M) \chi(Q)= (t-1)(t^2-4t+1)(t^4-1)$ and 
the largest real root  $2+\sqrt{3}$ 
is the stretch factor of $[f]$. 
\end{example}

\subsection{Construction of train track representatives}

Suppose that $[f]$ is pseudo-Anosov and 
$\mathfrak{g}: G \rightarrow G$ is a train track map of $[f]$. 
We define the {\it derivative} 
$D_{\mathfrak{g}}: \mathcal{E}_{\mathrm{dir}}(G) \rightarrow \mathcal{E}_{\mathrm{dir}}(G)$ as follows. 
Given an edge $ e \in \mathcal{E}_{\mathrm{dir}}(G)$, 
we write the edge path 
$\mathfrak{g}(e)= e_1 e_2 \dots e_{\ell}$, 
where $e_i \in \mathcal{E}_{\mathrm{dir}}(G)$. 
We set $D_{\mathfrak{g}}(e):= e_1$, that is, the initial edge of 
$\mathfrak{g}(e)$.

\begin{figure}[htbp]
\begin{center}
\includegraphics[width=4.5in]{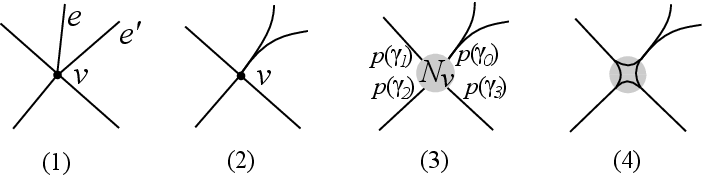}
\caption{Edges  $e$ and $e'$ are equivalent. 
Gray disks illustrate $N_v$. 
(1) Link $\mathrm{Lk}(v)$. 
(2) Step 1. 
(3) Step 2. 
(4) Step 3.}
\label{fig_gate}
\end{center}
\end{figure}

\begin{definition}
\label{def_gate} 
Consider the link 
$\mathrm{Lk}(v) \subset \mathcal{E}_{\mathrm{dir}}(G)$. 
Given $e_1, e_2 \in \mathrm{Lk}(v)$, we say that 
$e_1$ and $e_2$ are {\it equivalent} 
if $D_{\mathfrak{g}}^k(e_1) = D_{\mathfrak{g}}^k(e_2)$ for some $k \ge 1$. 
The equivalence classes are called {\it gates} 
incident to the vertex $v$. 
\end{definition}

\begin{definition}
\label{def_construction-track}
From a train track map 
$\mathfrak{g}: G \rightarrow G$ of $[f]$, 
we construct a train track 
$\tau= \tau(G)$ taking the following 3 steps 
at each vertex $v$ of $G$. 
(See Figure \ref{fig_gate}.)
\begin{enumerate} 
 \item[Step 1.] 
Deform the equivalent edges of $\mathrm{Lk}(v)$ in a neighborhood of  $v$ so that they are tangent at $v$. 

 \item[Step 2.] 
 Replace the vertex $v$ with a small disk $N_v$. 
 For each gate $\gamma$ incident to $v$, we assign a point $p(\gamma)$ on the boundary of  $N_v$. 

 \item[Step 3.] 
 If there exist an edge $e \in \mathcal{E}_{\mathrm{dir}}(G)$ and an integer 
 $k \ge 1$ such that 
 $\mathfrak{g}^k(e) = \cdots \overline{e_j} e_{\ell} \cdots $ for some 
 $e_j, e_{\ell} \in \mathrm{Lk}(v)$ with 
 $\gamma_j= [e_j] \ne \gamma_{\ell} = [e_ {\ell}]$, 
 then we join 
  $p(\gamma_j)$ and $p(\gamma_{\ell})$ by a smooth arc in $N_v$. 
 The arc intersects the boundary of $N_v$ transversely  
 at $p(\gamma_j)$ and $p(\gamma_{\ell})$, 
 and no two such arcs intersect in the interior of $N_v$. 
The resulting graph is $\tau=\tau(G)$. 
\end{enumerate}
\end{definition}

\begin{figure}[htbp]
\begin{center}
\includegraphics[width=4.3in]{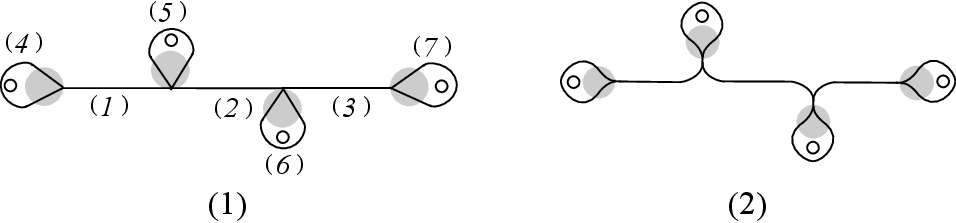}
\caption{
(1) Bestvina-Handel graph $G$ of Example \ref{example_4braid_BH}. 
(2) Train track $\tau(G)$.}
\label{fig_simple4braid}
\end{center}
\end{figure}

\begin{remark}
\label{remark_convension}
When we illustrate a Bestvina-Handel graph $G$, we may highlight the vertices with gray disks to show the additional structures of infinitesimal edges.
With this convention, $G$ can also represent 
the  track $\tau(G)$.  
For example, if $\mathfrak{g}: G \rightarrow G$ is 
the train track map of 
$[f]= \Gamma(\sigma_1 \sigma_3 \sigma_2^{-1})$ in Example \ref{example_4braid_BH}, 
Figure~\ref{fig_simple4braid}(2) depicts its train track $ \tau(G)$. 
It indicates that each vertex of $G$ is a partial vertex with three gates, and the graph in Figure~\ref{fig_simple4braid}(1) 
also represents the train track $\tau(G)$.
\end{remark}

We denote by $\mathcal{I}(\tau)$, 
the set of edges contained in $N_v$ 
for all vertices $v$ of $G$.  
Let $\mathcal{R}(\tau)$ be the set of other edges of $\tau$.

\begin{remark}
An edge in $\mathcal{R}(\tau)$ is called ``real" in \cite{BBK12}. 
This terminology is not the same as that 
in Definition \ref{definition_infinitesimal-real}. 
See also Remark \ref{remark_inconsistent}. 
\end{remark}

The set of all edges $\mathcal{E}(\tau)$ of $\tau$ is 
a union of $\mathcal{R}(\tau)$ and $\mathcal{I}(\tau)$: 
\begin{equation}
  \label{branch of tau}
  \mathcal{E}(\tau)=\mathcal{R}(\tau)\cup \mathcal{I}(\tau).
  \end{equation}
There is a natural surjection 
$$\pi: \tau \rightarrow G $$
that collapses all infinitesimal edges contained in $N_v$  
to the vertices $v$ of $G$ and that takes 
non-infinitesimal edges in $\mathcal{R}(\tau)$ 
to the corresponding edges of $G$. 
Since $\tau$ is $f$-invariant, 
the train track map $\mathfrak{g}: G \rightarrow G$ 
extends to a train  track representative 
$\mathfrak{t}: \tau \rightarrow \tau$ 
of $[f]$ with the following commutative diagram. 
$$
\begin{tikzcd}
    \tau \arrow[r, "\pi "] 
    \arrow[d, "\mathfrak{t}"]& G
    \arrow[d,"\mathfrak{g}"]\\
    \tau \arrow[r, "\pi "] & G
 \end{tikzcd}
$$

\begin{figure}[htbp]
\begin{center}
\includegraphics[width=4.5in]{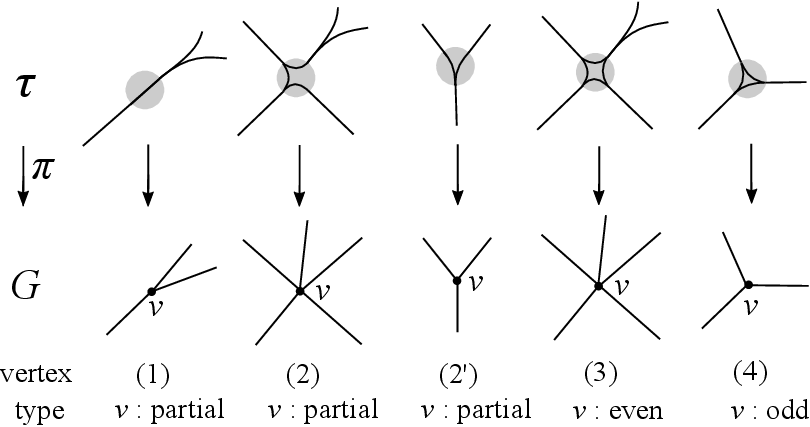}
\caption{
The surjection 
$\pi: \tau \rightarrow G$. 
(1) Case V1. (2)(2') Case V2. (3)(4) Case V3.}
\label{fig_vertex-type}
\end{center}
\end{figure}

It is proved in \cite[Proposition 3.3.4]{BestvinaHandel95} that 
if $v$ is a vertex of $G$ with $d$ gates, 
one of the following cases occurs. 
\begin{enumerate}
\item[Case V1.] 
$d=2$, and there exists a unique infinitesimal edge in $N_v$ that joins 
the two gates. 

\item[Case V2.]
$d>2$, and there exists an infinitesimal edge  that joins a pair of adjacent gates except for one pair. 
Infinitesimal edges in $N_v$ form a $d$-gon missing one side 
which is called  a {\it partial polygon} 
(Figure \ref{fig_vertex-type}(2)(2')).

\item[Case V3.]
$d >2$, and there exists an infinitesimal edge that joins 
a pair of adjacent gates. 
Infinitesimal edges in $N_v$ form a $d$-gon 
which is called an {\it infinitesimal polygon} 
(Figure \ref{fig_vertex-type}(3)(4)). 
\end{enumerate}

\begin{definition}
\label{definition_vertex-type}
Let $G$ be a Bestvina-Handel graph. 
\begin{enumerate}
\item 
Let $v$ be a vertex of $G$ with $d$ gates. 
We say that $v$ is {\it even} (resp. {\it odd}) 
if case V3 occurs and $d$ is even (resp. odd). 
The vertex $v$ is {\it partial} if case V1 or V2 occurs. 
See Figure \ref{fig_vertex-type}.  

\item 
We denote by $\mathcal{V}_{\mathrm{e}}(G)$, 
the set of all even vertices  of $G$, 
possibly $\mathcal{V}_{\mathrm{e}}(G) = \emptyset$. 
\end{enumerate}
\end{definition}

\begin{remark}
\label{remark_inconsistent}
Each edge in $\mathcal{I}(\tau)$ is infinitesimal 
(see Definition \ref{definition_infinitesimal-real}). 
On the other hand, 
an edge in $ \mathcal{R}(\tau)$ is not necessarily a real edge. 
In fact, the train track $\tau$ in Figure \ref{fig_simple4braid}(2) 
contains a $\mathfrak{t}$-invariant (peripheral) subgraph 
$P$ whose components are circles, 
representing the peripheral curves in $D_4$. 
Each edge in $P$ is an element of $\mathcal{R}(\tau)$ 
and  an infinitesimal edge. 
\end{remark}

Recall that 
the inclusion $G \rightarrow \Sigma$ is a homotopy equivalence. 
If $\mathfrak{g}: G\rightarrow G$ is a train track map, 
one can construct a pseudo-Anosov representative of the corresponding 
mapping class and the invariant foliations. 
See \cite[Section 3.4]{BestvinaHandel95}. 
The following lemma is a consequence of the Bestvina-Handel algorithm.

\begin{lemma}
\label{lemma_one-to-one} 
Let $f: \Sigma \rightarrow \Sigma$ be a pseudo-Anosov map 
with unstable and stable foliations $\mathcal{F}^u$ and $\mathcal{F}^s$. 
Let 
$\mathfrak{g}: G \rightarrow G$ be  a train track map of $[f]$. 
Then there is a one to one correspondence 
between the set of even (resp. odd) vertices of $G$ 
and the set of unpunctured even- (resp. odd-)pronged singularities of 
$\mathcal{F}^u$ and $\mathcal{F}^s$. 
\end{lemma}

\subsection{The weight spaces of Bestvina-Handel graphs}

Let $G$ be a Bestvina-Handel graph of $[f]$. 
Let $V(G) \simeq {\Bbb R}^{|\mathcal{E}(G)|}$ 
denote the vector space of real weights on the edges of $G$.  
The surjection $\pi_*\colon V(\tau) \to V(G)$ 
is induced by $\pi\colon \tau \to G$. 
Consider the weight space $W(\tau)$ of 
the train track $\tau$. 
The  \textit{weight space} $W(G)\subset V(G)$ of $G$ is defined by 
$$W(G):=\pi_*(W(\tau)).$$
The restriction $\pi_*|_{W(\tau)}: W(\tau) \rightarrow W(G)$ 
is denoted by 
$\pi_*: W(\tau) \rightarrow W(G)$.

\begin{figure}[htbp]
\begin{center}
\includegraphics[width=4.3in]{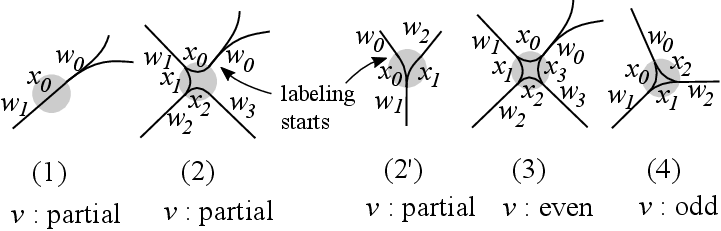}
\caption{
Gray disks illustrate $N_v$ for the vertex $v$. 
For $\eta\in V(G)$, 
the weight of the $i$th gate incident to  $v$ 
is denoted by 
$w_i:=w_i^{v}(\eta)$, and 
the weights on infinitesimal edges are denoted by $x_j$'s.}
\label{fig_weight-vertex-type}
\end{center}
\end{figure}

Let $v$ be a vertex of $G$ with $d$ gates. 
We label $0, \dots, d-1$ for  the gates counterclockwise. 
If $v$ is a partial vertex, 
the labeling starts at an end of the partial polygon 
as shown in Figure \ref{fig_weight-vertex-type}(2)(2'). 
Given an element $\eta\in V(G)$, 
let $w_i^{v}(\eta)$ be the weight (coefficient) 
of the $i$th gate incident to the vertex $v$, that is, 
$w_i^{v}(\eta)$ is the sum of weights (coefficients) 
of all edges in the $i$th gate.

\begin{lemma}[Lemma 2.3 in \cite{BBK12}]
\label{BBKLemma 2.3}
Let $\eta \in V(G)$. 
Then $\eta \in W(G)$ if and only if 
    for each non-odd vertex, that is, even  or partial vertex 
    $v$ with $d$ gates, 
    we have 
    $ \sum_{i=0}^{d-1}(-1)^i w_i^{v}(\eta)=0$. 
\end{lemma}

\begin{remark}
\label{rem_odd-vertex}
The odd vertices of $G$ are free in the statement of 
Lemma \ref{BBKLemma 2.3}. 
In fact, if $G$ contains an odd vertex $v$, 
a simple calculation shows that 
the switch condition is satisfied at every switch in 
$\tau(G)$ associated with the vertex $v$. 
(See also the proof of Lemma 2.3 in \cite{BBK12}.) 
For example, if $v$ is an odd vertex with three gates, 
then there is a $3$-pronged boundary component $\pi^{-1}(v)$ 
of $\tau(G)$. (See Definition \ref{definition_boundary-component}.)  
Given any real weights $w_0, w_1$ and $w_2$  at the three gates, 
there exist unique real weights $x_0, x_1$ and $x_2$ 
on the infinitesimal edges such that 
the switch conditions 
$ x_{i-1}+ x_i= w_i$ hold, 
where the index $i$ is considered $\mathrm{mod}\ 3$. 
See Figure \ref{fig_weight-vertex-type}(4). 
\end{remark}

\subsection{Homology polynomials}
\label{subsection_homology-polynomials}

Since the train track $\tau= \tau(G)$ is $f$-invariant, 
the linear map $\mathfrak{t}_{*}: V(\tau) \rightarrow V(\tau)$ 
preserves $W(\tau)$, that is,  
$\mathfrak{t}_{*}(W(\tau)) \subset W(\tau)$. 
We can define an induced linear map   
$\fg_*\colon W(G)\to W(G)$ by the following commutative diagram
$$
\begin{tikzcd}
    W(\tau) \arrow[r,"\pi_*"] \arrow[d,"\mathfrak{t}_*"] &W(G)\arrow[d,"\mathfrak{g}_*"] \\
     W(\tau) \arrow[r,"\pi_*"] &W(G)
\end{tikzcd}
$$

Birman-Brinkmann-Kawamuro stated 
 in Corollary 2.2 of \cite[p.31]{BBK12} 
that the characteristic polynomial 
$h(t):=\chi(\mathfrak{g}_*)$ of 
$\mathfrak{g}_*: W(G) \rightarrow W(G)$ 
is an invariant of the pseudo-Anosov mapping class $[f]$. 
Dodge-Kawamuro \cite{Kawamuro26} 
gave more details on the proof of this result. 

In \cite{BBK12}, $h(t)= \chi(\mathfrak{g}_*)$ is called 
the {\it homology polynomial} of $[f]$. 
We prove in Section \ref{subsection_invariance-homology} 
that the polynomial $h(t)$ is a conjugacy invariant of the pseudo-Anosov mapping class $[f]$. 
See also Remark \ref{remark_Dodge-Kawamuro}.

\begin{remark}
As stated in Theorem 1.1 of \cite{BBK12}, 
the homology polynomial $h(t)$ is a factor of $\chi(T)$, 
where $T$ is the transition matrix  of 
the train track map $\mathfrak{g}: G \rightarrow G$. 
To see this, 
we note that $\chi(T)$ is the characteristic polynomial of 
$\mathfrak{g}_*: V(G) \rightarrow V(G)$ 
with respect to the standard basis of $V(G)$. 
Since $W(G)$ is the subspace of $V(G)$ preserved by 
$\mathfrak{g}_*: V(G) \rightarrow V(G)$, 
the polynomial $h(t)(= \chi(\mathfrak{g_*}|_{W(G)}))$ 
should be a factor of $\chi(T)$. 
\end{remark}

\subsection{Puncture polynomials and symplectic polynomials}
\label{subsection_puncture-symplectic}

Following \cite{BBK12}, 
we introduce a skew-symmetric form on $W(G)$:  

\begin{definition}
\label{BBK Definition 3.2}
 For $\eta$, $\zeta \in W(G)$, we choose vectors $\eta'$, $\zeta'\in W(\tau)$ so that $\pi_*(\eta')=\eta$ and $\pi_*(\zeta')=\zeta$.
   A skew-symmetric form $\overline{\omega}$ on $W(G)$ is defined by
$$\overline{\omega}(\eta,\zeta):= \omega(\eta',\zeta'),
  $$ 
where $\omega$ is the Thurston symplectic form  on $W(\tau)$. 
See Section \ref{subsection_symplectic-form}. 
\end{definition}
The value of  $\overline{\omega}(\eta,\zeta)$ does not depend on the choices of $\eta', \zeta' \in W(\tau)$. 
To see this, it suffices to show that 
$\omega(\hat{\eta}, \hat{\zeta}) =0$ 
for any $\hat{\eta} \in \ker \pi_* \subset W(\tau)$ and 
$\hat{\zeta} \in W(\tau)$. 
This  computation is included in 
the proof of Proposition 3.1(1) in \cite{BBK12}.  
See also Remark \ref{remark_well-defined}.

 The \textit{radical} $Z(G)$ of $\overline{\omega}$ is defined by 
 $$Z(G):=\{ \eta \in W(G) \mid \overline{\omega}( \eta,\zeta)=0 
 \text{ for all }  \zeta \in W(G)  \}.$$
By definition, we have $\pi_*(Z(\tau))=Z(G)$.
Consider a direct sum decomposition 
$W(G) \simeq Z(G) \oplus W(G)/Z(G)$. 
Since  $\mathfrak{t}_*\colon W(\tau)\to W(\tau)$ preserves 
the Thurston symplectic form $\omega$ and 
$\mathfrak{g}_*$ is of the form 
$\mathfrak{g}_*=\pi_*\circ\mathfrak{t}_*\circ\pi_*^{-1},$  
the map $\mathfrak{g}_*: W(G) \rightarrow W(G)$ 
also preserves the form $\overline{\omega}$, 
and hence  $\mathfrak{g}_*(Z(G)) \subset Z(G)$. 
We can define the restriction map $\mathfrak{g}_*|_{Z(G)}$ 
and have the following commutative diagram.
$$
\begin{tikzcd}
    Z(\tau) \arrow[r,"\pi_*|_{Z(\tau)}"] \arrow[d,"\mathfrak{t}_*|_{Z(\tau)}"] &Z(G)\arrow[d,"\mathfrak{g}_*|_{Z(G)}"] \\
     Z(\tau) \arrow[r,"\pi_*|_{Z(\tau)}"] &Z(G)
\end{tikzcd}
$$
The quotient space $ W(G)/Z(G)$ inherits the form 
$\overline{\omega}$ 
and the induced map $\overline{\mathfrak{g}_*}$ 
preserves the form $\overline{\omega}$. 
The following commutative diagram holds.
\begin{equation*}
     \begin{tikzcd}
 0\arrow[r]& Z(G) \arrow[r] \arrow[d,"\mathfrak{g}_*|_{Z(G)}"]& W(G) \arrow[r,] \arrow[d,"\mathfrak{g}_*"]& W(G)/Z(G)\arrow[r] \arrow[d,"\overline{\mathfrak{g}_*}"]& 0\\
 0\arrow[r]& Z(G) \arrow[r] & W(G) \arrow[r,]& W(G)/Z(G)\arrow[r]& 0
\end{tikzcd}
\end{equation*}
The above diagram together with 
Lemma \ref{lemma_Rykken} yields the following.

\begin{proposition}[Theorem 3.1 in \cite{BBK12}]
\label{proposition_BBK-factorization}
The homology polynomial $h(t)= \chi(\mathfrak{g}_*)$ 
has a factorization 
$\chi(\mathfrak{g}_*)= \chi(\mathfrak{g}_*|_{Z(G)}) \chi(\overline{\mathfrak{g}_*})$. 
\end{proposition}

The polynomial 
$p(t):=  \chi(\mathfrak{g}_*|_{Z(G)})$ is called 
the \textit{puncture polynomial} and 
 $s(t):=  \chi(\overline{\mathfrak{g}_{*}})$ is called 
 the \textit{symplectic polynomial}. 
We collectively refer to $p(t)$, $s(t)$ and $h(t)=p(t)s(t)$ 
as the \textit{BBK polynomials} of $[f]$.

Birman-Brinkmann-Kawamuro stated  in Theorem 3.1 
of \cite[p.36]{BBK12} that 
the polynomials $p(t)$ and $s(t)$ are invariants of $[f]$, 
and the symplectic polynomial $s(t)$ 
contains the stretch factor of $[f]$ as its largest real root. 
In \cite{Kawamuro26}, 
more details on the proof of the invariance of 
the polynomials $p(t)$ and $s(t)$ are provided. 
We prove in Theorem \ref{theorem_coinv_puncturepoly} that 
the polynomials $p(t)$ and $s(t)$ are conjugacy invariants of $[f]$. 
See also Remark \ref{remark_Dodge-Kawamuro}. 

\section{Proof of conjugacy invariance of the BBK polynomials}
\label{section_conjugacy-invariance}

In this section, we introduce 
a  polynomial conjugacy invariant of pseudo-Anosov maps, which we call 
the unpunctured singularity polynomial $u(t)$. 
The polynomial $u(t)$ encodes the transitions   
among unpunctured even-pronged singularities of pseudo-Anosov maps. 
We prove that 
the BBK polynomials can be described by the train track polynomials 
together with $u(t)$. 
This gives us an alternative proof of the conjugacy invariance of 
the BBK polynomials.

Let $f: \Sigma \rightarrow \Sigma$ be a pseudo-Anosov map, 
where $\Sigma= \Sigma_{g,s}$ with $s \ge 1$. 
Let 
$\mathfrak{g}: G \rightarrow G$ be a train track map of $[f]$ 
and $\mathfrak{t}: \tau \rightarrow \tau$ the train track representative 
on $\tau= \tau(G)$. 
Several of the commutative diagrams that appear 
in  this section are based on the following commutative diagram 
(\ref{diagram_all}). 
The reader may find it helpful to keep this diagram 
(\ref{diagram_all}) 
in mind while reading this section.

\begin{equation}
\label{diagram_all}
\begin{tikzcd}
  & 0 \arrow[d]
  & 0 \arrow[d]
  & 0 \arrow[d] \\
0 \arrow[r]
  & \ker \pi_*|_{Z(\tau)} \arrow[r] \arrow[d, "\iota"]
  & Z(\tau) \arrow[r,"\pi_*|_{Z(\tau)}"]
  \arrow[d, "\iota"]
  & Z(G) \arrow[r]
      \arrow[d, "\iota"]
  & 0 \\
0 \arrow[r]
  & \ker \pi_* \arrow[r]
     \arrow[d]
  & W(\tau) \arrow[r, "\pi_*"]
      \arrow[d]
  & W(G) \arrow[r]
      \arrow[d]
  & 0 \\
0 \arrow[r]
  & \ker \overline{\pi_*} \arrow[r]
     \arrow[d]
  & W(\tau)/Z(\tau) \arrow[r, "\overline{\pi_*}"]
      \arrow[d]
  & W(G)/Z(G) \arrow[r]
      \arrow[d]
  & 0 \\
  & 0
  & 0
  & 0
\end{tikzcd}
\end{equation}
where both rows and columns are exact 
and $\iota$ denotes the inclusion. 
(We shall see 
$\ker \pi_*|_{Z(\tau)}= \ker \pi_*$, 
and hence 
$ \ker \overline{\pi_*} = \ker \pi_*/ \ker \pi_*|_{Z(\tau)}= \{0\}$.)

\subsection{Unpunctured singularity polynomials} 
\label{subsection_unpunctured}

We begin with the lemma.

\begin{lemma}
\label{lemma_kernel-pi}
The surjections 
$\pi_*: W(\tau) \rightarrow W(G)$ and 
$\pi_*|_{Z(\tau)}: Z(\tau) \rightarrow Z(G)$ 
satisfy the following. 
\begin{enumerate}
\item 
$\ker\pi_* = \ker \pi_*|_{Z(\tau)} $. 

\item 
$\ker \pi_*|_{Z(\tau)} $ has a basis 
$\{r_{\pi^{-1}(v)}\ |\  v\ \mathrm{is\  an\  even\  
vertex\  of\  }G\}$, 
where $r_{\pi^{-1}(v)}$ is the radical element 
constructed  in Definition \ref{definition_rc}. 
\end{enumerate}
\end{lemma}

\begin{proof}
Take an element 
$\zeta \in \ker \pi_* \subset W(\tau)$. 
By the definition of $\pi: \tau \rightarrow G$,  
the coefficient of 
$\zeta \in \ker \pi_*$ on any edge in $\mathcal{R}(\tau)$ is $0$. 
(For the definition of the set $\mathcal{R}(\tau)$, 
see the paragraph just after Remark \ref{remark_convension}.)  
Suppose that $e$ is an (infinitesimal) edge 
on  the infinitesimal polygon 
$\pi^{-1}(v)$ for some vertex $v$ of $G$. 
In case $v$ is an odd vertex, then $\zeta(e)=0$ 
(see Remark \ref{rem_odd-vertex}). 
In case  $v$ is a partial vertex, 
it is easy to see that $\zeta(e)=0$. 
Consider the remaining case where $v$ is an even vertex. 
Then $\pi^{-1}(v)$ is an even-pronged boundary component of $\tau$. 
Consider the radical element $r_{\pi^{-1}(v)}\in W(\tau)$. 
By the definition of $\pi: \tau \rightarrow G$, 
we have $r_{\pi^{-1}(v)} \in \ker \pi_*$. 
On the other hand, if there is a vector $\zeta\in \ker\pi_*$ where $\zeta(e)=a $ 
for an  edge $e \subset \pi^{-1}(v)$, 
then the coefficients of other edges $e' \subset \pi^{-1}(v)$ 
in $\zeta\in \ker\pi_*$ should be $a$ or $-a$ 
by the switch conditions, and it follows that 
$\zeta= \sum a_v r_{\pi^{-1}(v)}$ 
where $v$ rages over all even-pronged vertex of $G$. 
Hence, we have 
\begin{equation}
\label{equation_kernel}
 \ker\pi_*  
= \mathrm{span}\{r_{\pi^{-1}(v)}\ |\  v\ \mathrm{is\  an\  even\  
vertex\  of\  }G\}.   
\end{equation}
Since $r_{\pi^{-1}(v)} \in Z(\tau)$ 
(Proposition \ref{prop_prop8.1}), 
it holds $\ker\pi_* = \ker \pi_*|_{Z(\tau)} $. 
This completes the proof of (1). 

We turn to the proof of (2). 
Consider all even vertices $v$ of $G$. 
Then $\{r_{\pi^{-1}(v)}\}$ are linearly independent, 
since the infinitesimal polygons of the form $\pi^{-1}(v)$ 
are disjoint. 
This together with (\ref{equation_kernel}) implies claim (2). 
\end{proof}

\begin{remark}
\label{remark_well-defined}
The fact that 
the skew-symmetric form $\overline{\omega}$ on $W(G)$ 
is well defined 
follows from Lemma \ref{lemma_kernel-pi}. 
In fact, Lemma \ref{lemma_kernel-pi} tells us that 
if $\eta' \in \ker \pi_*$, then $\eta' \in Z(\tau)$. 
Hence, $\omega(\eta',\zeta')=0$ for any $\zeta' \in W(\tau)$. 
This implies that $\overline{\omega}$ is well defined. 
\end{remark}

By Lemma \ref{lemma_kernel-pi}, 
the dimension of 
$\ker \pi_* (= \ker \pi_*|_{Z(\tau)})$ 
equals $|\mathcal{V}_{\mathrm{e}}(G)|$, that is, 
the number of even vertices of $G$. 
Then the surjections 
$\pi_*: W(\tau) \rightarrow W(G)$ and 
$\pi_*|_{Z(\tau)}: Z(\tau) \rightarrow Z(G)$ yield 
the following result.

\begin{corollary}
\label{corollary_dimZ} 
Let $f: \Sigma \rightarrow \Sigma$ be a pseudo-Anosov map. 
Let $G$ be a Bestvina-Handel graph of $[f]$ 
and $\tau= \tau(G)$ the train track. 
We have the following 
about the dimensions of vector spaces. 
\begin{enumerate}
\item 
$\dim W(\tau) = \dim W(G) +|\mathcal{V}_{\mathrm{e}}(G)|$. 

\item 
$\dim Z(\tau) = \dim Z(G) +|\mathcal{V}_{\mathrm{e}}(G)|$.
\end{enumerate}
In particular, if $f$ is a fully punctured pseudo-Anosov map, 
then 
$\dim W(\tau) = \dim W(G)$ and 
$\dim Z(\tau) = \dim Z(G)$. 
\end{corollary}

Since $\mathfrak{g}$ preserves the set of even vertices, 
we have $\mathfrak{t}_ * (\ker \pi_*) \subset \ker \pi_*$. 
Hence, the restriction map 
$\mathfrak{t}_*|_{\ker \pi_*}: \ker \pi_* 
\rightarrow \ker \pi_*$ 
is defined.

\begin{theorem}
\label{theorem_kerpi*_coninv}
Let $f: \Sigma \rightarrow \Sigma$ be a pseudo-Anosov map, 
where $\Sigma= \Sigma_{g,s}$ with $s \ge 1$. 
The characteristic polynomial 
$\chi(\mathfrak{t}_*|_{\ker \pi_*})$ 
of $\mathfrak{t}_*|_{\ker \pi_*}: \ker \pi_* 
\rightarrow \ker \pi_*$ is a conjugacy invariant of 
$[f] \in \mathrm{MCG}(\Sigma)$. 
\end{theorem}

\begin{proof}
By Lemma \ref{lemma_kernel-pi}(1) 
we have 
$$
\chi(\mathfrak{t}_*|_{\ker \pi_*}) =  
\chi(\mathfrak{t}_*|_{\ker \pi_*|_{Z(\tau)}}), 
$$
 where 
$\mathfrak{t}_*|_{\ker \pi_*|_{Z(\tau)}}: 
\ker \pi_*|_{Z(\tau)} \rightarrow \ker \pi_*|_{Z(\tau)}$. 
We prove that 
$\chi(\mathfrak{t}_*|_{\ker \pi_*|_{Z(\tau)}})$ is a conjugacy invariant 
of $[f]$. 

Suppose that $\mathfrak{g}: G \rightarrow G$ and 
$\mathfrak{g}': G' \rightarrow G'$ are train track maps 
on Bestvina-Handel graphs $G$ and $G'$ for $[f]$, respectively. 
We have $f$-invariant train tracks $ \tau= \tau(G)$ 
and $ \tau'= \tau(G')$.  
Consider the train track representatives 
$ \mathfrak{t}: \tau \rightarrow \tau$ and 
$ \mathfrak{t}': \tau' \rightarrow \tau'$. 
We take the third train track $\tau''$ 
together with the train track representative 
$\mathfrak{t}'': \tau'' \rightarrow \tau''$ 
as in the proof of Theorem \ref{theorem_inv_t_*}.

Recall that the homotopy equivalence graph map 
$m: \tau \rightarrow \tau''$ 
yields the isomorphism 
$m_*|_{Z(\tau)}: Z(\tau) \rightarrow Z(\tau'')$ 
with the commutative diagram (\ref{equation_com_Z}). 
By Lemma \ref{lemma_kernel-pi} we have 
$$\ker \pi_*|_{Z(\tau)} 
= \mathrm{span}\{r_{\pi^{-1}(v)}\ |\ v\ \mbox{is an even vertex of }G\},$$ 
where $\pi: \tau \rightarrow G$ is the surjection. 
We take any even vertex $v$ of $G$ and consider 
the even-pronged boundary component $\pi^{-1}(v)$ of $\tau$. 
Propositions \ref{prop_prop5.3} and \ref{prop_prop8.1}(1) tell us that 
there is a boundary component $c''(v)$ of $\tau''$ 
corresponding to $\pi^{-1}(v)$ such that 
the radical element 
$r_{\pi^{-1}(v)} \in Z(\tau)$ maps to 
the radical element 
$r_{c''(v)} \in Z(\tau'')$ under the isomorphism 
$m_*|_{Z(\tau)}: Z(\tau) \rightarrow Z(\tau'')$.

We turn to the other homotopy equivalence graph map 
$m': \tau' \rightarrow \tau''$ 
which gives the isomorphism 
$m'_*|_{Z(\tau')}: Z(\tau') \rightarrow Z(\tau'')$. 
By Lemma  \ref{lemma_kernel-pi} again, 
we have 
$$\ker \pi'_*|_{Z(\tau')} 
= \mathrm{span}\{r_{\pi'^{-1}(v)}\ |\ v'\ \mbox{is an even vertex of }G'\},$$ 
where $\pi': \tau' \rightarrow G'$ is the surjection. 
By Lemma \ref{lemma_one-to-one}, 
there is a natural one to one correspondence 
between the set of even vertices of $G$ and 
the set of even vertices of $G'$ 
via the set of unpunctured even-pronged singularities of $\mathcal{F}^u$ 
and $\mathcal{F}^s$. 
For each even vertex $v$ of $G$, 
we pick an even vertex $v'$ of $G'$ 
corresponding to $v$, and consider 
the even-pronged boundary component $\pi'^{-1}(v')$ of $\tau'$. 
By Propositions \ref{prop_prop5.3}, 
the isomorphism 
$m'_*|_{Z(\tau')}: Z(\tau') \rightarrow Z(\tau'')$ 
takes  $r_{\pi'^{-1}(v')}$ to $r_{c''(v)} \in Z(\tau'')$ 
for the above boundary component $c''(v)$ of $\tau''$.  

Putting them together, one sees that 
for the above even vertices $v$ of $G$ and $v'$ of $G'$, 
the isomorphism 
$(m_*'|_{Z(\tau')})^{-1} \circ m_*|_{Z(\tau)}: 
Z(\tau) \rightarrow Z(\tau')$ 
takes the radical element $r_{\pi^{-1}(v)}$ to the radical element 
$r_{\pi'^{-1}(v')}$, 
and  gives the isomorphism 
between $\ker \pi_*|_{Z(\tau)}$ to $ \ker \pi'_*|_{Z(\tau')}$. 
Hence, the following  commutative diagram holds. 
$$
\begin{tikzcd}[column sep=2cm]
\ker \pi_*|_{Z(\tau)} \arrow[r,""] \arrow[d,"\mathfrak{t}_*|_{\ker \pi_*|_{Z(\tau)}}"']
  & \ker \pi'_*|_{Z(\tau')} \arrow[d,"
 \mathfrak{t}'_*|_{\ker \pi'_*|_{Z(\tau')}}"] \\
\ker \pi_*|_{Z(\tau)} \arrow[r,""]
  & \ker \pi'_*|_{Z(\tau')}
\end{tikzcd}
$$
Thus, 
 $\chi(\mathfrak{t}_*|_{\ker \pi_*|_{Z(\tau)}}) = 
 \chi(\mathfrak{t}'_*|_{\ker \pi_*|_{Z(\tau')}})  $ 
 is an invariant of $[f]$. 

The conjugacy invariance can be proved in the same manner 
as in the proof of Theorem \ref{theorem_inv_t_*}. 
This completes the proof. 
\end{proof}

We call $u(t):= \chi(\mathfrak{t}_*|_{\ker \pi_*})$ 
the {\it unpunctured singularity polynomial} of $[f]$.

\begin{corollary}
\label{corollary_ut1}
Let $f: \Sigma \rightarrow \Sigma$ be a pseudo-Anosov map, 
where $\Sigma= \Sigma_{g,s}$ with $s \ge 1$. 
Suppose that the unstable foliation $\mathcal{F}^u$ 
has no unpunctured even-pronged singularities.
Then  $u(t)=1$. 
In particular, if $f$ is fully punctured, then 
$u(t)=1$. 
\end{corollary}

\begin{proof}
By Lemmas \ref{lemma_one-to-one} and \ref{lemma_kernel-pi}, 
it follows that 
$\dim (\ker \pi_*)=0$. 
Hence, $u(t)= \chi(\mathfrak{t}_*|_{\ker \pi_*})=1$.  
\end{proof}

\subsection{Proof of conjugacy invariance of homology polynomials}
\label{subsection_invariance-homology}

Let $\mathfrak{g}\colon G \to G$ 
be a train track map of $[f]$. 
Let $\tau=\tau(G)$, and 
$\mathfrak{t}\colon \tau \to \tau$ 
a train track representative of $[f]$. 
We recall the linear maps 
$\mathfrak{g}_*\colon W(G) \to W(G)$ and 
$\mathfrak{t}_*\colon W(\tau) \to W(\tau)$. 
The following commutative diagram holds. 
$$\begin{tikzcd}
 0\arrow[r]& \ker \pi_* \arrow[r] \arrow[d,"\mathfrak{t}_*|_{\ker \pi_*}"]& W(\tau) \arrow[r, "\pi_*"] \arrow[d,"\mathfrak{t}_*"]& W(G)\arrow[r] \arrow[d,"\mathfrak{g}_*"]& 0\\
 0\arrow[r]& \ker \pi_* \arrow[r] & W(\tau) \arrow[r, "\pi_*"]& W(G)\arrow[r]& 0
\end{tikzcd}
$$
Applying Lemma \ref{lemma_Rykken} to the above diagram, 
we  obtain the following. 

\begin{proposition}\label{proposition_t and g}
The characteristic polynomial of 
$\mathfrak{t}_*: W(\tau) \rightarrow W(\tau)$ has a factorization 
$\chi(\mathfrak{t}_*)=
\chi(\mathfrak{g}_*)\chi(\mathfrak{t}_*|_{\ker \pi_*})$. 
\end{proposition}

\begin{theorem}
\label{theorem_coinv_homologypoly}
Let $f: \Sigma \rightarrow \Sigma$ be a pseudo-Anosov map, 
where $\Sigma= \Sigma_{g,s}$ with $s \ge 1$. 
The homology polynomial $h(t)= \chi(\mathfrak{g}_{*})$ 
is a conjugacy invariant of 
    $[f] \in \mathrm{MCG}(\Sigma)$. 
\end{theorem}
\begin{proof}
    By Proposition \ref{proposition_t and g}, we can write $\chi(\mathfrak{g}_{*})$ as $\chi(\mathfrak{g}_{*})=
    \chi(\mathfrak{t}_*)/\chi(\mathfrak{t}_*|_{\ker \pi_*})$. 
By Theorems \ref{theorem_inv_t_*} and \ref{theorem_kerpi*_coninv}, 
the polynomials $\chi(\mathfrak{t}_*)$ and 
$\chi(\mathfrak{t}_*|_{\ker \pi_*})$ 
are  conjugacy invariants of $[f]$. 
Thus, $\chi(\mathfrak{g}_{*})$ is  a conjugacy invariant of $[f]$.
\end{proof}

\subsection{Proof of conjugacy invariance of puncture/symplectic polynomials}
\label{subsection_invariance-puncture}

Recall that 
$\overline{\omega}$ is the skew-symmetric form on $W(G)$, 
and $Z(G)$ is the radical of $\overline{\omega}$. 
Since $\mathfrak{t}_*\colon W(\tau)\to W(\tau)$ 
preserves the Thurston symplectic form $\omega$ on $W(\tau)$, we have the following commutative diagram. 
$$
     \begin{tikzcd}
 0\arrow[r]& \ker \pi_*|_{Z(\tau)} \arrow[r] \arrow[d,"\mathfrak{t}_*|_{\ker \pi_*|_{Z(\tau)}}"]& Z(\tau) \arrow[r, "\pi_*|_{Z(\tau)}"] \arrow[d,"\mathfrak{t}_*|_{Z(\tau)}"]& Z(G)\arrow[r] \arrow[d,"\mathfrak{g}_*|_{Z(G)}"]& 0\\
 0\arrow[r]& \ker \pi_*|_{Z(\tau)} \arrow[r] & Z(\tau) \arrow[r, "\pi_*|_{Z(\tau)}"]& Z(G)\arrow[r]& 0
\end{tikzcd}
$$
Recall that 
$\chi(\mathfrak{t}_*|_{\ker \pi_*}) =  
\chi(\mathfrak{t}_*|_{\ker \pi_*|_{Z(\tau)}})$. 
Applying Lemma \ref{lemma_Rykken} to the above diagram, 
we obtain the following.

\begin{proposition}\label{proposition_t and gZ}
The characteristic polynomial $\chi(\mathfrak{t}_*|_{Z(\tau)})$ 
of $\mathfrak{t}_*|_{Z(\tau)}: Z(\tau) \rightarrow Z(\tau)$ 
has a factorization 
$\chi(\mathfrak{t}_*|_{Z(\tau)})=
\chi(\mathfrak{g}_*|_{Z(G)})\chi(\mathfrak{t}_*|_{\ker \pi_*})$. 
\end{proposition}

Let $\overline{\mathfrak{t}_*}\colon W(\tau)/Z(\tau)\to W(\tau)/Z(\tau)$ be the map induced by $\mathfrak{t}_*: W(\tau) \rightarrow W(G)$, 
and let 
$\overline{\mathfrak{g}_*}\colon W(G)/Z(G)\to W(G)/Z(G)$ 
be the map induced by $\mathfrak{g}_*: W(G) \rightarrow W(G)$. 
We have the following commutative diagram.
\begin{equation}\label{com_dia_W/Z}
     \begin{tikzcd}
 0\arrow[r]& \ker \overline{\pi_*} \arrow[r] \arrow[d,"\overline{\mathfrak{t}_*}|_{\ker \overline{\pi_*}}"]& W(\tau)/Z(\tau) \arrow[r, "\overline{\pi_*}"] \arrow[d,"\overline{\mathfrak{t}_{*}}"]& W(G)/Z(G)\arrow[r] \arrow[d,"\overline{\mathfrak{g}_{*}}"]& 0\\
 0\arrow[r]& \ker \overline{\pi_*} \arrow[r] & W(\tau)/Z(\tau) \arrow[r, "\overline{\pi_*}"]& W(G)/Z(G)\arrow[r]& 0
\end{tikzcd}
\end{equation}
By Lemma \ref{lemma_kernel-pi}(1),  
$\ker \overline{\pi_{*}} = \ker \pi_*/ \ker \pi_*|_{Z(\tau)} = \{0\}$. 
See also  diagram (\ref{diagram_all}). 
Hence, the map 
$\overline{\pi_*}: W(\tau)/Z(\tau) \rightarrow W(G)/Z(G) $ 
in (\ref{com_dia_W/Z}) 
(induced by $\pi_*: W(\tau) \rightarrow W(G)$) 
is an isomorphism, and we obtain the following.

\begin{proposition}
\label{proposition_symmplectic-polynomial}
$\chi(\overline{\mathfrak{t}_{*}})= \chi(\overline{\mathfrak{g}_{*}}) 
(= s(t))$. 
\end{proposition}

\begin{theorem}
\label{theorem_coinv_puncturepoly}
Let $f: \Sigma \rightarrow \Sigma$ be a pseudo-Anosov map, 
where $\Sigma= \Sigma_{g,s}$ with $s \ge 1$. 
The puncture polynomial $p(t)= \chi(\mathfrak{g}_{*}|_{Z(G)})$ 
and the symplectic polynomial $s(t)= \chi(\overline{\mathfrak{g}_*})$ 
are conjugacy invariants of 
    $[f] \in \mathrm{MCG}(\Sigma)$.
\end{theorem}

\begin{proof}
    By Proposition \ref{proposition_t and gZ}, 
    we can write 
    $p(t)= \chi(\mathfrak{g}_{*}|_{Z(G)})$ as $$\chi(\mathfrak{g}_{*}|_{Z(G)})=
    \chi(\mathfrak{t}_*|_{Z(\tau)})/\chi(\mathfrak{t}_*|_{\ker \pi_*}).$$
Then $p(t)$ is a conjugacy invariant, since  
the radical polynomial $\chi(\mathfrak{t}_*|_{Z(\tau)})$ 
and the unpunctured singularity polynomial 
$\chi(\mathfrak{t}_*|_{\ker \pi_*})$ 
are conjugacy invariants of $[f]$ 
(see Theorems \ref{theorem_invariant-tt-puncture} and \ref{theorem_kerpi*_coninv}). 

The conjugacy invariance of $s(t)= \chi(\overline{\mathfrak{g}_*})$ 
follows from Theorem \ref{theorem_invariant-tt-puncture} and 
Proposition \ref{proposition_symmplectic-polynomial}. 
\end{proof}

Theorem \ref{theorem_characterization} is a straightforward modification 
of the following result, which is a characterization of 
the BBK polynomials. 

\begin{theorem}
\label{proposition_poly_tau}
Let $f: \Sigma \rightarrow \Sigma$ be a pseudo-Anosov map, 
where $\Sigma= \Sigma_{g,s}$ with $s \ge 1$. 
Let $\mathfrak{g}\colon G\to G$ be a train track map 
and $\mathfrak{t}\colon \tau \to \tau$ 
the train track representative of $[f]$. 
Then we have the following. 
\begin{enumerate}
\item
$h(t)(= \chi(\mathfrak{g}_*)) = \chi(\mathfrak{t}_*)/ \chi(\mathfrak{t}_*|_{\ker\pi_*})$

\item 
$p(t) (= \chi(\mathfrak{g}_*|_{Z(G)})) = 
\chi(\mathfrak{t}_*|_{Z(\tau)})/ \chi(\mathfrak{t}_*|_{\ker\pi_*})$

\item 
$s(t) (= \chi(\overline{\mathfrak{g}_*}))= \chi(\overline{\mathfrak{t}_*})$.
\end{enumerate}
In particular, if $f$ is fully punctured, 
then 
$h(t) = \chi(\mathfrak{t}_*)$, 
$p(t) = \chi(\mathfrak{t}_*|_{Z(\tau)})$ 
and 
$s(t) = \chi(\overline{\mathfrak{t}_*})$. 
\end{theorem}

\begin{proof}
For  claims (1), (2) and (3), 
see Propositions \ref{proposition_t and g}, 
\ref{proposition_t and gZ} and 
\ref{proposition_symmplectic-polynomial}. 
If $f$ is fully punctured, then 
$u(t) (=\chi(\mathfrak{t}_*|_{\ker\pi_*}))=1$ 
(Corollary \ref{corollary_ut1}). 
This completes the proof.
\end{proof}

\begin{proof}[Proof of Corollary \ref{corollary_invariance}]
The conjugacy invariance of the BBK polynomials follows from 
Theorems \ref{theorem_coinv_homologypoly} and 
\ref{theorem_coinv_puncturepoly}. 
\end{proof}

\begin{definition}
\label{definition_reciprocal}
Let $Q(t)= \sum_{i=1}^n a_i t^{\alpha_i}$ be a polynomial. 
Its {\it degree} and its {\it reciprocal polynomial} are defined by 
$ \mathrm{deg}(Q):= \max \alpha_i$ and 
$Q_*(t):= t^{\mathrm{deg}(Q)} Q(t^{-1})$. 
We say that 
$Q(t)$ is {\it reciprocal} if $Q(t)= \pm Q_*(t)$; 
it is {\it palindromic} if $Q(t)= Q_*(t)$, 
and  {\it anti-palindromic} if  $Q(t) = -Q_*(t)$. 
\end{definition}

We now prove the proposition stated 
in Section \ref{section_introduction}. 

\begin{proof}[Proof of Proposition \ref{prop_reciprocal}]
Suppose that $f$ has at least two singularity orbits. 
In the proof of  Theorem 6.1 of \cite{HironakaTsang25}, 
it is proved that 
the polynomials 
$\chi(\overline{\mathfrak{t}_*})$, 
$\chi(\mathfrak{t}_*|_{Z(\tau)}) $ and 
$\chi(\mathfrak{t}_*)$ are reciprocal 
for a particular train track 
representative $\mathfrak{t}: \tau \rightarrow \tau$. 
Since the train track polynomials 
are conjugacy invariants of pseudo-Anosov mapping classes $[f]$ 
(Theorems \ref{theorem_inv_t_*} and \ref{theorem_invariant-tt-puncture}), 
the same properties hold for all train track representatives 
of $[f]$. 
In the case where $f$ has a single singularity orbit, 
one can prove the statement in the same manner.

As an alternative proof, it is proved in \cite{BBK12} that 
the BBK polynomials 
$s(t)= \chi(\overline{\mathfrak{g}_*})$, 
$p(t)= \chi(\mathfrak{g}_*|_{Z(G)}) $ and 
$h(t)=\chi(\mathfrak{g}_*)$ are reciprocal. 
If $f$ is a fully punctured pseudo-Anosov map, 
the BBK polynomials coincide with the train track polynomials 
(Theorem \ref{proposition_poly_tau}). 
Therefore, the train track polynomials 
$\chi(\overline{\mathfrak{t}_*})$, 
$\chi(\mathfrak{t}_*|_{Z(\tau)}) $ and 
$\chi(\mathfrak{t}_*)$ are reciprocal.  
\end{proof}

The following lemma can be found in \cite[Lemma 2.5]{BBK12}.

\begin{lemma}
\label{BBKLemma 2.5}
The dimension of $W(G)$ is 
$$
\begin{array}{ll}
\#\{\text{edges of }G \}-\#\{\text{vertices of }G \}+1
& \mbox{if}\ \tau \  \mbox{is orientable},  
\\
\#\{\text{edges of }G \}-\#\{\text{non-odd vertices of }G\}
& \mbox{if}\ \tau\ \mbox{is non-orientable}. 
\end{array}
$$
\end{lemma}

\begin{figure}[htbp]
\begin{center}
\includegraphics[width=3.8in]{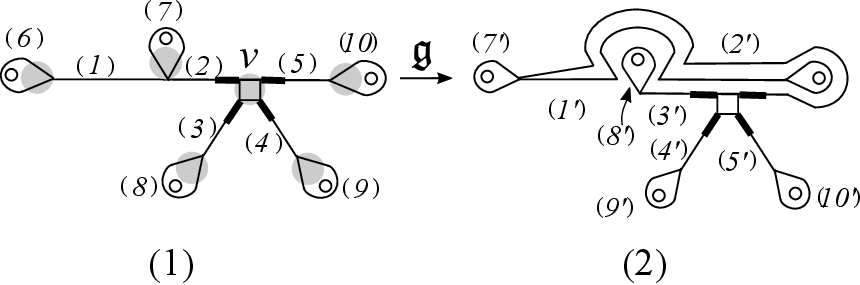}
\caption{
Example \ref{ex_running}. 
Bold edges are half edges emanating from the even vertex $v$.}
\label{fig_g_1-3}
\end{center}
\end{figure}

\begin{figure}[htbp]
\begin{center}
\includegraphics[width=3.8in]{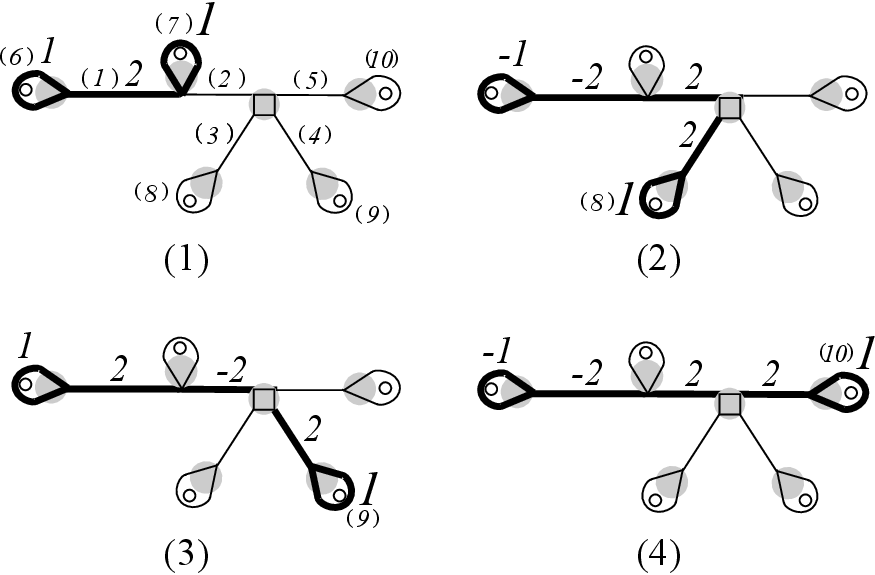}
\caption{Vectors (1) $\eta_{(7)}$, (2) $\eta_{(8)}$, (3) $\eta_{(9)}$ 
and (4) $\eta_{(10)}$ in $W(G)$. 
An edge without number means its weight is $0$.}
\label{fig_basis_g_1-3}
\end{center}
\end{figure}

\begin{example}
\label{ex_running}
Consider  
$[f]= \Gamma(\sigma_1 \sigma_2^{-1} \sigma_3^{-1} \sigma_4^{-1}) 
\in \mathrm{MCG}(D_5)$ 
which is pseudo-Anosov. 
We compute the BBK polynomials. 
Figure \ref{fig_g_1-3} illustrates 
a train track map $\mathfrak{g}: G \rightarrow G$ for $[f]$.  
The edges of the  graph $G$ are labeled 
$(1), \dots, (10)$, and 
the set of vertices consists of  an even vertex $ v$  and 
five partial vertices.   
See the gray disks in Figure \ref{fig_g_1-3}(1) 
and Remark \ref{remark_convension}. 

Since  $\tau(G)$ is non-orientable, 
Lemma \ref{BBKLemma 2.5} tells us  $\dim W(G)=10-6=4$. 
Take vectors $\eta_{(i)} \in V(G)$ 
$(i= 7,8,9,10)$ 
as shown in Figure \ref{fig_basis_g_1-3}. 
For example, 
$\eta_{(7)}= \left(\begin{smallmatrix}
2 & 0 & 0 & 0 & 0 & 1 & 1 & 0 & 0 & 0
\end{smallmatrix}\right)$. 
By Lemma \ref{BBKLemma 2.3}, 
one sees that 
$\eta_{(i)} \in W(G)$. 
Among them, the vector $\eta_{(i)}$ is the only one 
such that the edge labeled  $(i)$ has non-zero weight (coefficient). 
Hence, they are linearly independent and 
yield a basis of $ W(G)$. 
The representation matrix  of 
$\mathfrak{g}_*: W(G) \rightarrow W(G)$ 
with respect to this basis is given by 
$ \left(\begin{smallmatrix}
2 & 1 & 0 &0
\\ 0 & 0& 1 & 0 
\\ 0 & 0 & 0& 1
\\ 1& 1& -1& 1
\end{smallmatrix}\right)$. 
Hence, 
$ h(t)=\chi(\mathfrak{g}_*)=t^4-3t^3+3t^2-3t+1$. 
Let $\mathfrak{t}: \tau \rightarrow \tau$ be the train track 
representative of $[f]$ obtained from 
$\mathfrak{g}: G \rightarrow G$. 
All the boundary components of $\tau$ have $1$-pronged, 
except the infinitesimal polygon $\pi^{-1}(v)$ 
which is $4$-pronged. 
From the transition of half edges emanating from $v$ 
(Figure \ref{fig_g_1-3}), 
we see that 
$\mathfrak{t}_*(\gamma_{\pi^{-1}(v)}) = -\gamma_{\pi^{-1}(v)}$. 
(See Definition \ref{definition_rc} for the radical element.)
Hence, 
the radical polynomial $\chi(\mathfrak{t}|_{Z(\tau)})$ and 
the unpunctured singularity polynomial 
$u(t)=\chi(\mathfrak{t}_*|_{\ker \pi_*})$ 
satisfy 
$\chi(\mathfrak{t}|_{Z(\tau)}) = u(t)= t+1$. 
This together with Proposition \ref{proposition_t and gZ} 
implies that  $p(t)=1$. 
By factorization 
$h(t)= p(t)s(t)$ (Proposition \ref{proposition_BBK-factorization}), 
we obtain 
$s(t)(= h(t))= t^4-3t^3+3t^2-3t+1$. 
The largest real root  equals the stretch factor 
$\lambda \approx 2.15372$ of $[f]$. 
\end{example}

\section{Examples}
\label{section_examples}

Let $f: \Sigma \rightarrow \Sigma$ be a pseudo-Anosov map. 
We define a pseudo-Anosov map $f_{\mathrm{even}}$ obtained from $f$. 
To do this, let $\mathcal{S}_{\mathrm{e}}$ 
be the set of all unpunctured even-pronged singularities of $\mathcal{F}^u$, possibly $\mathcal{S}_{\mathrm{e}}= \emptyset$. 
We set $\Sigma_{\mathrm{e}}:= \Sigma \setminus \mathcal{S}_{\mathrm{e}}$. 
Then we have a pseudo-Anosov map 
$$f_{\mathrm{even}}:= f|_{\Sigma_{\mathrm{e}}}: 
\Sigma_{\mathrm{e}} \rightarrow \Sigma_{\mathrm{e}}$$ 
that is obtained from $f$ by removing all unpunctured 
even-pronged singularities. 
If there are no unpunctured odd-pronged singularities, then 
$f_{\mathrm{even}}= f^{\circ}$ is a fully punctured pseudo-Anosov map. 

\begin{theorem}
\label{theorem_unpunctured-even}
Let $f: \Sigma \rightarrow \Sigma$ be a pseudo-Anosov map, 
where $\Sigma= \Sigma_{g,s}$ with $s \ge 1$. 
Then the pseudo-Anosov maps 
$f_{\mathrm{even}}: \Sigma_{\mathrm{e}} \rightarrow \Sigma_{\mathrm{e}}$ 
and $f^{\circ}: \Sigma^{\circ} \rightarrow \Sigma^{\circ}$ 
have the same BBK polynomials. 
More precisely, the following properties hold. 
\begin{enumerate}
\item 
The symplectic polynomials of $[f_{\mathrm{even}}]$, $[f^{\circ}]$ 
and $[f]$ are the same. 

\item 
The puncture polynomials of $[f_{\mathrm{even}}]$ and $[f^{\circ}]$ 
are both equal to $p(t)u(t)$, where $p(t)$ and $u(t)$ are the puncture polynomial and the unpunctured singularity polynomial of $[f]$, respectively.
\end{enumerate}
\end{theorem}

\begin{proof}
The claim (1) holds from Propositions 
\ref{prop_same-ttp} and \ref{proposition_symmplectic-polynomial}. 

We turn to the proof of (2) for $[f_{\mathrm{even}}]$. 
The proof of (2) for $[f^{\circ}]$ is similar. 
By Proposition \ref{proposition_t and gZ},  
the radical polynomial 
$\chi(\mathfrak{t}_*|_{Z(\tau)})$ of $[f]$ is given by 
$p(t)u(t) 
(= \chi(\mathfrak{g}_*|_{Z(G)})\chi(\mathfrak{t}_*|_{\ker \pi_*}))$. 
Let $p'(t)$ be the puncture polynomial of $[f_{\mathrm{even}}]$. 
By Corollary \ref{corollary_ut1}, 
the unpunctured singularity polynomial of $[f_{\mathrm{even}}]$ is $1$. 
This together with Proposition \ref{proposition_t and gZ} 
implies that the radical polynomial of $[f_{\mathrm{even}}]$ is 
equal to the puncture polynomial $p'(t)$ 
of $[f_{\mathrm{even}}]$. 
Since the radical polynomials of  $[f_{\mathrm{even}}]$ 
and $[f]$ are the same (Proposition \ref{prop_same-ttp}), 
we obtain $p'(t)= p(t) u(t)$. 
This completes the proof.
\end{proof}

\begin{example}
Let $f: D_5 \rightarrow D_5$ be a pseudo-Anosov map 
that represents 
$\Gamma(\sigma_1 \sigma_2^{-1} \sigma_3^{-1} \sigma_4^{-1}) 
\in \mathrm{MCG}(D_5)$.  
See Example \ref{ex_running}. 
Consider the fully punctured pseudo-Anosov map 
$f^{\circ}: D_6 \rightarrow D_6$. 
Let $p^{\circ}(t)$ and $s^{\circ}(t)$ be the puncture polynomial 
and the symplectic polynomial of  $[f^{\circ}]$, respectively. 
By Theorem \ref{theorem_unpunctured-even} 
together with Example \ref{ex_running}, 
we have 
$p^{\circ}(t)= p(t) u(t) = t+1$ and 
$s^{\circ}(t)= t^4-3t^3+3t^2-3t+1$. 
\end{example}

The next example tells us that 
puncture polynomials are  useful, 
since one can detect that mapping classes are not conjugate.

\begin{figure}[htbp]
\begin{center}
\includegraphics[width=4.5in]{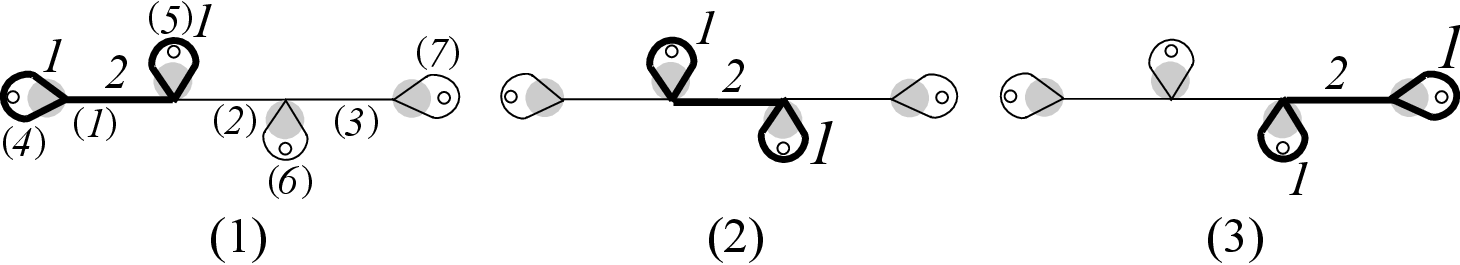}
\caption{Vectors (1) $\eta_{(1)}$, (2) $\eta_{(2)}$ and (3) $\eta_{(3)}$ 
 in $W(G)$.}
\label{fig_another_basis4braid}
\end{center}
\end{figure}

\begin{example}
\label{example_puncture-polynomials}
Let $\Delta_4 = \sigma_1 \sigma_2 \sigma_3 
\sigma_1 \sigma_2 \sigma_1 \in B_4$ be the positive half twist. 
Consider the two pseudo-Anosov mapping classes 
$\phi:= \Gamma(\sigma_1 \sigma_3 \sigma_2^{-1}) $ 
and $\phi' = \Gamma(\sigma_1 \sigma_3 \sigma_2^{-1} \Delta_4)$ 
in $\mathrm{MCG}(D_4)$. 
The former mapping class $\phi$ is the one 
in Example \ref{example_4braid_BH}. 
We claim that the puncture polynomials 
can detect $\phi$ and $\phi'$ are not conjugate, 
although their symplectic polynomials are the same: 
The graph $G$ in Figure \ref{fig_g_1-1-1}(1) is 
a Bestvina-Handel graph of $\phi$, 
 and it is also a Bestvina-Handel graph of $\phi'$. 
 Note that $G$ is preserved by the $\pi$ rotation 
 $r: G \rightarrow G$ about the middle point on the edge 
 labeled (2). 

Recall the train track map 
$\mathfrak{g}: G \rightarrow G$ of $[\phi]$ 
as in Figure \ref{fig_g_1-1-1}. 
One can check that the graph map 
$\mathfrak{g}':=r \circ \mathfrak{g}: G \rightarrow G$ 
is a train track map of $[\phi']$. 
Let $\mathfrak{t}: \tau \rightarrow \tau$ 
and $\mathfrak{t}': \tau \rightarrow \tau$ 
be the train track representatives of $\phi$ and $\phi'$, respectively. 
All the boundary components of the common train track 
$\tau= \tau(G)$ are $1$-pronged, 
except for the $2$-boundary component $c$ 
(Figure \ref{fig_radical-element}(2)). 
Consider the radical element $r_c$ 
(Figure \ref{fig_radical-element}(3)). 
By computations, one sees that 
$\mathfrak{t}_*(r_c)= r_c$ and 
$\mathfrak{t}'_*(r_c)= -r_c$. 
This tells us that 
the radical polynomials of $\phi$ and $\phi'$ are 
$t-1$ and $t+1$, respectively. 
Since $\phi$ (resp. $\phi'$) is 
a fully punctured pseudo-Anosov mapping class, 
by Theorem \ref{proposition_poly_tau} 
the radical polynomial is equal to the puncture polynomial 
of $\phi$ (resp. $\phi'$). 
Hence, the puncture polynomials of $\phi$ and $\phi'$ differ, 
and this means that 
$\phi, \phi \in \mathrm{MCG}(D_4)$ are not conjugate.   

Let us compute the symplectic polynomials. 
Since $\tau$ is non-orientable, 
by Lemma \ref{BBKLemma 2.5} we have 
 $\dim W(G)=7-4=3$.  
The vectors $\eta_{(i)} $ 
$(i= 1,2,3)$  
as shown in Figure \ref{fig_another_basis4braid} 
is a basis of $ W(G)$.
If $M$ and $M'$ are the representation matrices 
of $\mathfrak{g}_*: W(G) \rightarrow W(G)$ 
and $\mathfrak{g}'_*: W(G) \rightarrow W(G)$, respectively, 
then 
$$M = 
\left(\begin{smallmatrix}
1 & 1 & 0 
\\ 1 & 3 & 1 
\\ 0 & 1 & 1
\end{smallmatrix}\right)\  \mbox{and} \ 
M' = 
\left(\begin{smallmatrix}
0 & 0 & 1 
\\ 0 & 1 & 0 
\\ 1 & 0 & 0
\end{smallmatrix}\right)
\left(\begin{smallmatrix}
1 & 1 & 0 
\\ 1 & 3 & 1 
\\ 0 & 1 & 1
\end{smallmatrix}\right)
=
\left(\begin{smallmatrix}
0 & 1 & 1 
\\ 1 & 3 & 1 
\\ 1 & 1 & 0
\end{smallmatrix}\right)$$
where $M$ is the block matrix of the transition matrix $T$ 
of $\mathfrak{g}: G \rightarrow G$ 
as in Example \ref{example_4braid_BH}. 
Hence, the homology polynomials $\phi$ and $\phi'$ are 
$\chi(M)=(t-1)(t^2-4t+1)$ and 
$\chi(M')=(t+1)(t^2-4t+1)$, respectively. 
By the factorization of the homology polynomial 
(Proposition \ref{proposition_BBK-factorization}),  
$\phi$ and $\phi'$ have the same symplectic polynomial 
$t^2-4t+1$. 
\end{example}

\begin{figure}[htbp]
\begin{center}
\includegraphics[width=4in]{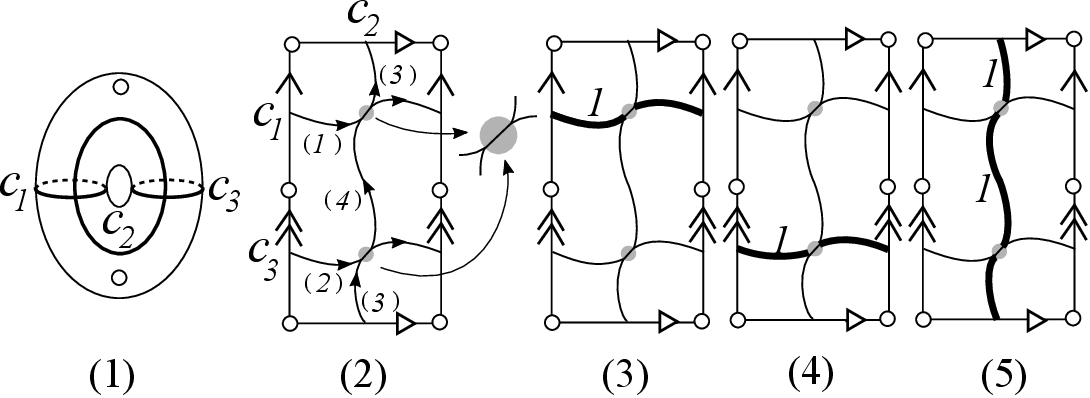}
\caption{
(1) Curves $c_i$. (2) Orientable train track $\tau= \tau(G)$. 
Vectors (3) $\eta_{(1)}$, (4) $\eta_{(2)}$ and (5) $\eta_{(3)}$ 
in $W(\tau)$. 
An edge without number means its weight is $0$.}
\label{fig_twice_track}
\end{center}
\end{figure}

\begin{example}
\label{example_twice-punctured-torus}

On the twice-punctured torus $\Sigma_{1,2}$, 
let $\delta_i$ be the right-handed Dehn twist 
about the simple closed curve $c_i$ $(i=1,2,3)$ 
as in Figure \ref{fig_twice_track}(1). 
Consider a rectangle with its four corners and 
the two middle points of its sides removed, 
and identify the sides of the rectangle according to the arrows 
of the same kinds 
as shown in Figure \ref{fig_twice_track}(2). 
We then obtain a twice-punctured torus $\Sigma_{1,2}$. 
Consider 
$[f]= \delta_2^{-1} \circ \delta_3 \circ \delta_1 \in \mathrm{MCG}(\Sigma_{1,2})$. 
This is a pseudo-Anosov mapping class obtained from Penner's construction \cite{Penner88}. 
Figure \ref{fig_twice_track}(2) illustrates 
an invariant train track $\tau$ which is orientable. 
Collapsing the two gray disks into points, we obtain a 
Bestvina-Handel graph $G$ of $[f]$. 
By Lemma \ref{BBKLemma 2.5}, 
$\dim W(G)= 4-2+1=3$. 
Since the pseudo-Anosov map 
$f: \Sigma_{1,2} \rightarrow \Sigma_{1,2}$ is fully punctured, 
by Corollary \ref{corollary_dimZ} 
we have $ \dim W(\tau)= \dim W(G) =3$.  
Take  vectors $\eta_{(i)}$ $(i= 1,2,3)$ 
as shown in Figure \ref{fig_twice_track}. 
The set of vectors $\{\eta_{(i)}\}$ is a basis of $W(\tau)$. 
Let $\mathfrak{t}: \tau \rightarrow \tau$ be 
the train track representative of $[f]$. 
The representation matrix of 
$\mathfrak{t}_*: W(\tau) \rightarrow W(\tau)$ is 
$\left(\begin{smallmatrix}
1 & 1 & 0 
\\ 1 & 3 & 1 
\\ 0 & 1 & 1
\end{smallmatrix}\right) 
= 
\left(\begin{smallmatrix}
1 & 0 & 0 
\\ 1 & 1 & 1 
\\ 0 & 0 & 1
\end{smallmatrix}\right) 
\left(\begin{smallmatrix}
1 & 0 & 0 
\\ 0 & 1 & 0
\\ 0 & 1 & 1
\end{smallmatrix}\right)
\left(\begin{smallmatrix}
1 & 1 & 0 
\\ 0 & 1 & 0
\\ 0 & 0 & 1
\end{smallmatrix}\right)$. 
Hence, by Theorem \ref{proposition_poly_tau}, 
$\chi(\mathfrak{t}_*)= \chi(\mathfrak{g}_*) 
(=h(t))= (t-1)(t^2-4t+1)$. 
By Proposition \ref{prop_prop8.1}(3), 
$\dim Z(\tau)= 2-1=1$, and 
by the factorization of the homology polynomial 
(Proposition \ref{proposition_BBK-factorization}), 
we have 
$\chi(\mathfrak{t}_*|_{Z(\tau)}) = \chi(\mathfrak{g}_*|_{Z(G)}) 
(= p(t)) = t-1$ and 
$s(t)= t^2-4t+1$. 

\end{example}

\subsection*{Acknowledgments} 
We are grateful to Keiko Kawamuro for kindly explaining the results in 
\cite{BBK12}, which motivated the present work. 
We would like to thank Chi Cheuk Tsang for patiently answering our questions concerning the results in \cite{HironakaTsang25}. 
Their comments and explanations were very helpful to us.
We would also like to thank Shunsuke Kano and Merrick R. Dodge  for helpful discussions.

\subsection*{Funding} 
The first author  was supported by JSPS KAKENHI Grant Number 
JP21K03247, JP23K22396, JP23K25778, JP26K06807. 

\subsection*{Tool Use}
ChatGPT (OpenAI) and Google Gemini were used to check English grammar, improve the wording of the manuscript, and obtain information on LaTeX syntax and usage.

\bibliographystyle{alpha}
\bibliography{BBKinvariant}
\end{document}